\documentclass[12pt]{amsart}
\usepackage{amsmath,amssymb,latexsym,cite,mathrsfs}
\usepackage{verbatim,wasysym,cite}
\usepackage[left=2.8cm,right=2.8cm,top=2.9cm,bottom=2.9cm]{geometry}

\usepackage[utf8]{inputenc}
\usepackage{microtype}
\usepackage{color,enumitem,graphicx}
\usepackage[colorlinks=true,urlcolor=blue, citecolor=red,linkcolor=blue,
linktocpage,pdfpagelabels, bookmarksnumbered,bookmarksopen]{hyperref}

\usepackage[english]{babel}

\renewcommand{\epsilon}{{\varepsilon}}

\numberwithin{equation}{section}
\newtheorem{theorem}{Theorem}[section]
\newtheorem{lemma}[theorem]{Lemma}
\newtheorem{remark}[theorem]{Remark}
\newtheorem{definition}[theorem]{Definition}
\newtheorem{proposition}[theorem]{Proposition}
\newtheorem{corollary}[theorem]{Corollary}
\DeclareMathOperator{\eq}{eq}
\DeclareMathOperator{\rad}{rad}
\DeclareMathOperator{\dom}{dom}
\DeclareMathOperator{\re}{Re}
\DeclareMathOperator{\im}{Im}
\DeclareMathOperator{\supp}{supp}
\title[ Instability of ground states]
{Instability of ground states for the NLS equation with  potential on  the star graph}

\author[Alex H. Ardila]{}
\email{ardila@impa.br}
\author[Liliana Cely]{}
\email{mlcelyp@ime.usp.br}
\author[Nataliia Goloshchapova]{}
\email{nataliia@ime.usp.br}

\subjclass[2010]{Primary: 35Q55; Secondary: 35Q40}
\keywords{Nonlinear Schr\"{o}dinger equation, linear potential, generalized Kirchhoff's condition, ground state, orbital stability}

\begin{document}

\maketitle

\centerline{\scshape Alex H. Ardila}
{\footnotesize
\centerline{Universidade Federal de Minas Gerais}
\centerline{ CEP 30123-970, Belo Horizonte-MG, Brazil}
}
\centerline{\scshape Liliana Cely}
{\footnotesize
\centerline{ Universidade de S\~ao  Paulo}
\centerline{CEP 05508-090,  Cidade Universit\'aria,  S\~ao  Paulo-SP, Brazil}
}
\centerline{\scshape Nataliia Goloshchapova}
{\footnotesize
\centerline{ Universidade de S\~ao  Paulo}
\centerline{CEP 05508-090,  Cidade Universit\'aria, S\~ao  Paulo-SP, Brazil}
}

\begin{abstract}
We study the nonlinear Schr\"odinger equation with an arbitrary real  potential $V(x)\in (L^1+L^\infty)(\Gamma)$  on a star graph $\Gamma$.  At the vertex an interaction occurs described by the generalized Kirchhoff condition with strength $-\gamma<0$. We show the existence of ground states $\varphi_{\omega}(x)$ as minimizers of the action functional on the Nehari manifold under additional negativity and decay conditions on $V(x)$. Moreover, for $V(x)=-\dfrac{\beta}{x^\alpha}$, in the supercritical case, we prove that the standing waves $e^{i\omega t}\varphi_{\omega}(x)$ are orbitally unstable in $H^{1}(\Gamma)$ when  $\omega$ is large enough. Analogous result holds for an arbitrary $\gamma\in\mathbb{R}$ when  the standing waves have  symmetric profile.

\end{abstract}

\medskip

\section{Introduction}
We consider the following focusing nonlinear Schr\"odinger equation on an infinite star graph  $\Gamma$:
\begin{equation}\label{NLS}
\begin{cases} 
i\partial_{t}u(t,x)=-\Delta_{\gamma}u(t,x)+V(x)u(t,x)-\left|u(t,x)\right|^{p-1}u(t,x),\quad(t,x)\in \mathbb{R}\times\Gamma,\\
u(0,x)=u_{0}(x),
\end{cases} 
\end{equation}
where $\gamma>0$, $p>1$, $u:\mathbb{R}\times\Gamma\rightarrow\mathbb{C}^{N}$,  and  $\Delta_{\gamma}$ is the Laplace operator with the generalized Kirchhoff condition at the vertex of  $\Gamma$ ($\cdot'$ stands for spatial derivative): $$v_1(0)=\ldots =v_N(0),\quad \sum^{N}_{e=1}v^{\prime}_{e}(0)=-\gamma v_{1}(0).$$ We assume that the potential $V(x)=(V_e(x))_{e=1}^N$ is real-valued and   satisfies the \textit{Assumptions} (see notation section):
\begin{itemize}
    \item[{\bf 1.}]\,\textit{Self-adjointness assumption}: $V(x)\in L^1(\Gamma)+L^\infty(\Gamma).$
\item[{\bf 2.}]\,\textit{Weak continuity assumption}: $\lim\limits_{x\to\infty}V_e(x)=0.$
\item[{\bf 3.}]\,\textit{Minimizing  assumption}: $\int\limits_{\mathbb R^+}V_e(x)|\phi(x)|^2dx<0$ for all $\phi(x)\in H^1(\mathbb{R}^+) \setminus\{0\}$.
\item[{\bf 4.}]\,\textit{Virial identity assumption}: $xV'(x)\in L^1(\Gamma)+L^\infty(\Gamma)$. 
  \end{itemize}
 Notice that \textit{Assumption 3} essentially guarantees  $(Vu,u)_2<0,\, u\in H^1(\Gamma)\setminus\{0\},$ and $V(x)\leq 0$  a.e. on $\Gamma$ (see Remark \ref{V_neg}).
 
NLS equation  \eqref{NLS} models wave propagation in thin waveguides (we refer the reader to \cite{BecImp15,QGAH,JolSem11, Kuc02} for the details).
 The study of stability properties of the multi-dimensional NLS with a linear potential
 \begin{align*}
     & i\partial_{t}u(t,x)=-\Delta u(t,x)+V(x)u(t,x)-\left|u(t,x)\right|^{p-1}u(t,x),\\&(t,x)\in \mathbb{R}\times\mathbb{R}^n,\quad 1+4/n\leq p< 1+4/(n-2),\end{align*}
 was initiated  in \cite{RosWei88}. More precisely, the authors proved orbital stability of $e^{i\omega t}\varphi_\omega(x)$ for $\omega$  sufficiently close to minus the smallest eigenvalue of the operator $-\Delta+V$ (under the assumptions $V(x)\in L^\infty(\mathbb{R}^n)$, $\lim\limits_{|x|\to \infty}V(x)=0$.)  In \cite{STW}, the stability results obtained by \cite{RosWei88} were improved for  $V(x)$ satisfying more general assumptions. 
 
 Recently in \cite{Oht18}, the author studied strong instability (by blow-up) of the standing waves in the case of harmonic potential $V(x)=|x|^2.$  In particular, he proved  strong instability under certain concavity condition for the associated  action functional (cf. Theorem \ref{INLS} below). The same idea was applied in \cite{FukOht19} to investigate strong instability for $V(x)=-\dfrac{\beta}{|x|^\alpha}, \, 0<\alpha<\min\{2,n\},\,\beta>0.$ The reader is also referred to \cite{Mar14} for more information about NLS near soliton dynamics.  
 
In the case $V(x)\equiv 0$,  the well-posedness in $H^{1}(\Gamma)$, variational and stability/instability properties  of  \eqref{NLS} have been extensively studied during the last decade. The well-posedness  results were obtained in  \cite{AQFF,  NataMasa2020}, whereas the  existence,  stability and variational properties  of ground states were studied in \cite{AQFF, NQTG, AnGo2018, AnGo20182, Kai19}. Moreover, the regularity and strong instability  results were  elaborated in \cite{NataMasa2020}. 

On the other hand, the NLS with potential on graphs is  little studied. To our  knowledge, the only results concerning the existence and stability of standing waves  were obtained in \cite{CAFiNo2017, Cacci2017, Ardila2018}. 
  In the subcritical ($1<p<5$) and critical ($p=5$) case orbitally stable standing waves $e^{i\omega t}\varphi_{\omega}(x)$  were constructed in  \cite{CAFiNo2017, Cacci2017} under specific conditions on $V(x).$
  Subsequently in \cite{Ardila2018} the orbital stability of $e^{i\omega t}\varphi_{\omega}(x)$ was studied in the supercritical case ($p>5$). More precisely,  it was shown (by solving a local energy-minimization problem) that  $e^{i\omega t}\varphi_{\omega}(x)$  is  stable when the mass of $\varphi_{\omega}(x)$ is sufficiently small. 
  
  \vspace{0.4cm}
  
  In this paper, we show  the existence and orbital instability of the standing wave solutions to \eqref{NLS} relying on methods developed in \cite{FukOht19,STW_a}. Moreover, we state regularity of the solutions to the Cauchy problem for the initial data from the domain of the operator $-\Delta_\gamma+V(x).$ This result is used to show virial identity which  is the key ingredient in the proof of the  instability result.
\subsection{Notation}
 We consider a graph $\Gamma$ consisting of a central vertex $\nu$ and $N$ infinite half-lines attached to it.  One may identify $\Gamma$ with the disjoint union of the intervals $I_{e}=(0,\infty)$, $e=1,\ldots,N$, augmented by the central vertex $\nu=0$. Given a function $v:\Gamma\rightarrow\mathbb{C}^{N}$, $v=(v_{e})^{N}_{e=1}$, where $v_{e}:(0,\infty)\rightarrow\mathbb{C}$ denotes the restriction of $v$ to $I_{e}$. We  denote by $v_{e}(0)$ and $v'_{e}(0)$ the limits of $v_{e}(x)$ and  $v'_{e}(x)$  as $x\rightarrow 0^{+}$. 
 
 We say that a function $v$ is continuous on $\Gamma$ if every restriction $v_{e}$ is continuous on $I_{e}$ and $v_{1}(0)=\ldots=v_{N}(0)$. The space of continuous functions is denoted by $C(\Gamma).$ 

The natural Hilbert space associated to  the Laplace operator $\Delta_{\gamma}$  is $L^{2}(\Gamma)$, which is
defined as $L^{2}(\Gamma)=\bigoplus^{N}_{e=1}L^{2}(\mathbb{R}^{+})$, and is equipped with the norm 
\begin{align*}
\left\|v\right\|^2_2=\int\limits_{\Gamma}\left|v\right|^{2}\,dx=\sum^{N}_{e=1}\int\limits_{0}^\infty\left|v_{e}(x)\right|^{2}\,dx.
\end{align*}
The inner product in $L^2(\Gamma)$ is denoted by $(\cdot,\cdot)_2.$
The space $L^{q}(\Gamma)$ for $1\leq q \leq\infty $ is defined analogously, and  $\|\cdot\|_q$ stands for its norm.
The Sobolev spaces $H^{1}(\Gamma)$ and $H^{2}(\Gamma)$ are  defined as
\begin{align*}
H^{1}(\Gamma)=\left\{v\in C(\Gamma):\,\,v_{e}\in H^{1}(\mathbb{R}^{+}),\,\,\,e=1,\ldots,N\right\},\\
H^{2}(\Gamma)=\left\{v\in C(\Gamma):\,\,v_{e}\in H^{2}(\mathbb{R}^{+}),\,\,\,e=1,\ldots,N\right\}.
\end{align*}
We consider the self-adjoint operator $ H_{\gamma, V}$ on $L^{2}(\Gamma)$: 
\begin{equation}\label{Ku}
\begin{split}
&(H_{\gamma, V}v)_{e}=-(\Delta_{\gamma}v)_{e}+V_ev_{e}=-v''_{e}+V_ev_{e},\\
&\dom( H_{\gamma, V})=\left\{v\in H^{1}(\Gamma):\,\,\,-v''_{e}+V_ev_{e}\in L^{2}(\mathbb{R}^{+}),\,\,\,\sum^{N}_{e=1}v^{\prime}_{e}(0)=-\gamma v_{1}(0)\right\}.
\end{split}
\end{equation}
 When $\gamma=0$, the condition at the vertex in \eqref{Ku} is usually referred as free or Kirchhoff boundary condition. For $\gamma\in\mathbb{R}$ the operator $ H_{\gamma, V}$ has a precise interpretation as the self-adjoint operator on $L^{2}(\Gamma)$ uniquely associated with the closed semibounded quadratic form $ {F}_{\gamma, V}$ defined on $H^{1}(\Gamma)$ by (see Lemma \ref{Ape1} in Appendix)
\begin{equation}\label{F}
\begin{split}
 {F}_{\gamma, V}(v)&=\|v'\|^{2}_2-\gamma\left|v_{1}(0)\right|^{2}+(Vv,v)_2\\
&=\sum^{N}_{e=1}\int\limits_{0}^\infty\left|v'_{e}(x)\right|^{2}\,dx-\gamma\left|v_{1}(0)\right|^{2}+\sum^{N}_{e=1}\int\limits_{0}^\infty V_e(x)|v_e(x)|^2\,dx.
\end{split}
\end{equation}
Note that we can formally rewrite  \eqref{NLS} as
\begin{equation*}
i\partial_{t}u(t)=E^{\prime}(u(t)),
\end{equation*}
where $E$ is the energy functional defined by
\begin{align*}
E(u)=\frac{1}{2} {F}_{\gamma, V}(u)-\frac{1}{p+1}\|u\|^{p+1}_{p+1}.
\end{align*}
The energy functional is well-defined on $H^{1}(\Gamma)$ since the potential $V(x)$ belongs to $(L^{1}+L^{\infty})(\Gamma)$ (see  Lemma \ref{Ape1} in Appendix).
\subsection{Standing waves and instability results}
By a standing wave of \eqref{NLS}, we mean a solution of the form $e^{i\omega t}\varphi(x)$, where $\omega\in\mathbb{R}$ and $\varphi$ is a solution of the stationary equation
\begin{equation}\label{S}
 H_{\gamma, V}\phi+\omega\phi-\left|\phi\right|^{p-1}\phi=0.
\end{equation}
We define two functionals on $H^{1}(\Gamma)$:
\begin{align*}
S_{\omega}(v):&=\frac{1}{2} {F}_{\gamma, V}(v)+\frac{\omega}{2}\left\|v\right\|^{2}_2-\frac{1}{p+1}\|v\|^{p+1}_{p+1}\quad (\text{\textit{action functional}}),\\
I_{\omega}(v):&= {F}_{\gamma, V}(v)+\omega\left\|v\right\|^{2}_2-\|v\|^{p+1}_{p+1}.
\end{align*}
Observe that \eqref{S} is equivalent to $S^\prime_{\omega}(\phi)=0$ (see \cite[Theorem 4]{AQFF}) and $I_{\omega}(v)=\partial_{\lambda}S_{\omega}(\lambda v)\left.\right|_{_{\lambda=1}}=\langle S^\prime_{\omega}(v) ,v\rangle$. Denote the set of non-trivial solutions to \eqref{S} by 
\begin{equation*}
\mathcal{B}_{\omega}=\Big\{v\in  H^{1}(\Gamma)\backslash \{0\}\,:\,\,S^\prime_{\omega}(v)=0\Big\}.
\end{equation*}
A ground state for \eqref{S} is a function $\varphi\in\mathcal{B}_{\omega}$ that minimizes $S_{\omega}$ on $\mathcal{B}_{\omega}$, and the set of  ground states is given  by
\begin{equation*}
\mathcal{G}_{\omega}=\Big\{\phi\in\mathcal{B}_{\omega}\,:\,\,S_{\omega}(\phi)\leq S_{\omega}(v)\,\,\mbox{for all}\,\,v\in\mathcal{B}_{\omega}\Big\}.
\end{equation*}
We consider the minimization problem on the Nehari manifold
\begin{align*}\label{d}
d_{\omega}&=\inf\left\{S_{\omega}(v)\,:\,\,v\in H^{1}(\Gamma)\backslash \{0\},\,\,I_{\omega}(v)=0\right\},
\end{align*}
and the set of minimizers
\begin{align*}
\mathcal{M}_{\omega}&=\left\{\phi\in H^{1}(\Gamma)\backslash \{0\}\,:\,\,S_{\omega}(\phi)=d_{\omega},\,\,\,I_{\omega}(\phi)=0\right\}.
\end{align*}

We now state the first result, which provides the existence of the minimizer for $d_{\omega}$ when  the  strength $-\gamma$ is sufficiently strong.  Denote (see  Lemma \ref{lemma_4.11})
\begin{equation}\label{omega_0}
-\omega_0:=\inf\sigma( H_{\gamma, V})=\min\sigma_p(H_{\gamma, V}) <0.   
\end{equation}
\begin{proposition}\label{GS}
Let $p>1$,  $\omega>\omega_{0},$ and $V(x)=\overline{V(x)}$ satisfy Assumptions 1-3. Then there exists $\gamma^{\ast}>0$ such that the set  $\mathcal{G}_{\omega}$ is not empty for any $\gamma>\gamma^{\ast}$, in particular,  $\mathcal{G}_{\omega}=\mathcal{M}_{\omega}$. If $\varphi_{\omega}\in\mathcal{G}_{\omega}$, then there exist $\theta\in\mathbb{R}$ and a positive function $\phi\in \dom(H_{\gamma, V})$ such that $\varphi_{\omega}(x)=e^{i\theta}\phi(x)$.
\end{proposition}
\noindent To be precise, $\gamma^*$  is given in \cite{AQFF} by \begin{equation}\label{gamma_star}\int\limits_0^1(1-t^2)^{\frac{2}{p-1}}dt=\frac{N}{2}\int\limits_{\frac{\gamma^*}{N\sqrt{\omega}}}^1(1-t^2)^{\frac{2}{p-1}}dt.\end{equation}  
The condition $\gamma>\gamma^*$ guarantees that the action functional  $S_\omega$ constrained to the Nehari manifold admits an absolute minimum when $V(x)\equiv 0$.
\begin{remark}\label{V_neg} The proof of the last assertion of Proposition \ref{GS} essentially uses that $V(x)\leq 0$ a.e. on $\Gamma$, which is a consequence of   \textit{Assumption 3}. 
 
To show this one  observes that $\int\limits_{\mathbb R^+} -V_e(x)\phi(x)dx\geq 0$ for all nonnegative functions $\phi(x)$ from $C_c(\mathbb{R}^+)$ (the set of continuous functions with compact support). Indeed, let $\tilde\phi(x)$ be an extension  onto $\mathbb R$ by zero of a nonnegative function $\phi(x)\in C_c(\mathbb{R}^+)$. Take
 $\{\phi_n(x)\}\subset C_c^\infty(\mathbb R)$
  such that $\phi_n\underset{n\to\infty}{\longrightarrow}\sqrt{\tilde\phi}$ uniformly, and $\supp \tilde\phi, \supp\phi_n\subset K\subset\mathbb R_+$, where $K$ is a compact set. Then $\phi^2_n\underset{n\to\infty}{\longrightarrow}\tilde\phi$ uniformly, and, by the  Dominated Convergence Theorem, we  get $$-\int\limits_{\mathbb R^+} V_e(x)\phi_n^2(x)dx\underset{n\to\infty}{\longrightarrow}-\int\limits_{\mathbb R^+} V_e(x)\phi(x)dx\geq 0.$$
  Now, since  $f(\phi)=-\int\limits_{\mathbb R^+} V_e(x)\phi(x)dx$ is a  positive linear functional on $C_c(\mathbb R^+)$, then, by the Riesz–Markov–Kakutani  representation theorem for positive linear functionals, we conclude the  existence of a unique Radon measure  $\mu$ on $\mathbb R^+$ such that $f(\phi)=\int\limits_{\mathbb R^+} \phi(x)d\mu(x)$. On the other hand, $f(\phi)=\int\limits_{\mathbb R^+} v(x)\phi(x)d\nu(x),$ where $\nu(A)=\int\limits_{A} |V_e|dx$ for $A$ from the Borel $\sigma$-algebra on $\mathbb R^+$, and $v(x)=\left\{\begin{array}{c}
    \frac{V_e(x)}{|V_e(x)|},\,\, x\in\{x: V_e(x)\neq 0\}  \\
     0,\qquad \,\, \text{otherwise}.
\end{array}\right.$ Finally, from the uniqueness stated in \cite[Theorem 2.5.12]{Fed69}  it follows that $\mu=\nu$ and  $v=1$ $\nu$-a.e. on $\mathbb{R}^{+},$ hence $-V_{e} \geq 0$  $\nu$-a.e. on $\mathbb{R}^{+}$. This  implies $-V_e \geq 0$ Lebesgue-a.e. on $\mathbb{R}^+$ since  the Lebesgue measure and the measure $\nu$ are mutually absolutely continuous on the set $\{x:V_e(x)\neq 0\}$.
\end{remark}

The next step in the study of ground states for \eqref{S} is to investigate their stability properties. We define orbital stability  as follows.
\begin{definition}  \label{D1} For $\varphi_{\omega}\in\mathcal{G}_{\omega}$, we set
\begin{equation}  \label{N}
N_{\delta}(\varphi_{\omega}):=\big\{v\in H^{1}(\Gamma)\,:\,\,\inf_{\theta\in\mathbb{R}}\left\|v-e^{i\theta}\varphi_{\omega}\right\|_{H^{1}(\Gamma)}<\delta\big\}.
\end{equation}
We say that a standing wave solution $e^{i\omega t}\varphi_{\omega}(x)$ of \eqref{NLS} is orbitally stable in $H^{1}(\Gamma)$ if for any $\epsilon>0$ there exists $\delta>0$ such that for any $u_{0}\in N_{\delta}(\varphi_{\omega})$, the solution $u(t)$ of \eqref{NLS}   satisfies $u(t)\in N_{\epsilon}(\varphi_{\omega})$ for all $t\geq0$. Otherwise, $e^{i\omega t}\varphi_{\omega}(x)$ is said to be orbitally unstable in $H^{1}(\Gamma)$.
\end{definition}	
Using the ideas developed in \cite{STW_a,FukOht19}, we obtain  a sufficient condition for the instability of standing waves when $p>5$ (supercritical case). The main result of this paper is the following:
\begin{theorem} \label{INLS}
Assume that $p>5$,   $\gamma>\gamma^*$,  $\omega>\omega_0$, and  $V(x)=\overline{V(x)}$ satisfies  Assumptions 1-4.   If $\varphi_{\omega}(x)\in\mathcal{G}_{\omega}$ and  $\partial^{2}_{\lambda}E(\varphi^{\lambda}_{\omega})\left|_{\lambda=1}\right.<0$, where $\varphi^{\lambda}_{\omega}(x):=\lambda^{1/2}\varphi_{\omega}(\lambda x)$ for $\lambda>0$, then the standing wave solution $e^{i\omega t}\varphi_{\omega}(x)$ of \eqref{NLS} is orbitally unstable in $H^{1}(\Gamma)$.
\end{theorem}
To prove  Theorem \ref{INLS} we use the variational  characterization given in Proposition \ref{GS} and  virial identity \eqref{well_18}.  
Notice that  the standing wave solution $e^{i\omega t}\varphi_{\omega}(x)$ of \eqref{NLS} with $\gamma>0$ and $V(x)\equiv 0$ is unstable in $H^{1}(\Gamma)$ when $p>5$ and $\omega$ is large enough (see \cite[Remark 6.1]{AQFF} and also \cite[Theorem 1.4]{NataMasa2020}). Below  we state  that this also holds true  for $\gamma>0$ and slowly decaying potential $V(x)=\dfrac{-\beta}{x^\alpha},\, 0<\alpha<1,\, \beta>0$  (i.e. $\partial^{2}_{\lambda}E_{\omega}(\varphi^{\lambda}_{\omega})\left|_{\lambda=1}\right.<0$  for sufficiently large $\omega$). The choice of the potential is due to its ``homogeneity" property, which is principal for the proof (see formula \eqref{E1}).   
\begin{corollary}\label{CINLS}
Assume that    $V(x)=\dfrac{-\beta}{x^\alpha}$, $\beta>0$,  $0<\alpha<1$, $\gamma>\gamma^*$,  $p>5$. If $\varphi_{\omega}(x)\in\mathcal{G}_{\omega}$, then there exists $\omega^*=\omega^*(\beta,\alpha,\gamma, p)\in(\omega_0,\infty)$ such that for any $\omega\in (\omega^*,\infty)$ the standing wave solution $e^{i\omega t}\varphi_{\omega}(x)$ of \eqref{NLS} is orbitally unstable in $H^{1}(\Gamma)$.  
\end{corollary}
As far as we know, these are  the first  results on  instability of ground states for the NLS with potential on graphs.
In Subsection \ref{subsec_eq}, we state the counterparts to  Proposition \ref{GS}, Theorem \ref{INLS}, Corollary \ref{CINLS} in the space $H^1_{\eq}(\Gamma)$ of symmetric functions and arbitrary $\gamma\in\mathbb{R}.$

The  paper is organized as follows. In Section \ref{S:1}, we  prove Proposition \ref{PP1} that concerns local well-posedness in the energy domain. In Section \ref{S:2}, we provide  the proof of Proposition \ref{GS} .  Section \ref{S:3} is devoted to the  proof of  Theorem \ref{INLS} and Corollary \ref{CINLS}. In the Appendix we discuss some properties of the operator $ H_{\gamma, V}$.
\section{local existence results and virial identity}
\label{S:1}
We start with the  proof of  the following key lemma involving the estimate of $H^1$-norm of the unitary group generated by the self-adjoint operator $ H_{\gamma, V}.$ 
\begin{lemma}
Let $e^{-i H_{\gamma, V}t}$ be a unitary group generated by $ H_{\gamma, V}$. Then $e^{-i H_{\gamma, V}t}H^1(\Gamma)\subseteq H^1(\Gamma)$ and 
\begin{equation}\label{group}
  \|e^{-i H_{\gamma, V}t}v\|_{H^1(\Gamma)}\leq C \|v\|_{H^1(\Gamma)}.  
\end{equation}
\end{lemma}
\begin{proof}
The idea of the proof was given in \cite{CAFiNo2017}
 (see formula (2.5)). However, some additional technical details seem useful.

Let $m>\omega_0$, where $\omega_0$ is  given by \eqref{omega_0}.
Remark that $H^1(\Gamma)=\dom\left( {F}_{\gamma, V}\right)=\dom(( H_{\gamma, V}+m)^{1/2})$ (see, for instance, \cite[Chapter VI, Problem 2.25]{Kat95}). Since $e^{-i H_{\gamma, V}t}$ is bounded, we get  for $v\in H^1(\Gamma)$ $$e^{-i H_{\gamma, V}t}( H_{\gamma, V}+m)^{1/2}v=( H_{\gamma, V}+m)^{1/2}e^{-i H_{\gamma, V}t}v.$$ Hence $e^{-i H_{\gamma, V}t}v\in H^1(\Gamma)$ and  $e^{-i H_{\gamma, V}t}H^1(\Gamma)\subseteq H^1(\Gamma).$ Further, using $L^2$-unitarity of $e^{-i H_{\gamma, V}t}$, we obtain for $v\in H^1(\Gamma)$ 
\begin{align*}
 & {F}_{\gamma, V}(v)+m\|v\|_2^2=\left(( H_{\gamma, V}+m)^{1/2}v,( H_{\gamma, V}+m)^{1/2}v\right)_2\\&=\left(e^{-i H_{\gamma, V}t}( H_{\gamma, V}+m)^{1/2}v, e^{-i H_{\gamma, V}t}( H_{\gamma, V}+m)^{1/2}v\right)_2\\&=\left(( H_{\gamma, V}+m)^{1/2}e^{-i H_{\gamma, V}t}v, ( H_{\gamma, V}+m)^{1/2}e^{-i H_{\gamma, V}t}v\right)_2= {F}_{\gamma, V}(e^{-i H_{\gamma, V}t}v)+m\|e^{-i H_{\gamma, V}t}v\|_2^2. 
\end{align*}
From the proof of Lemma \ref{lemma_4.11}-$(ii)$ we get
\begin{align*}
 &C_2\|e^{-i H_{\gamma, V}t}v\|^2_{H^1(\Gamma)}\leq  {F}_{\gamma, V}(e^{-i H_{\gamma, V}t}v)+m\|e^{-i H_{\gamma, V}t}v\|_2^2\\& = {F}_{\gamma, V}(v)+m\|v\|_2^2\leq C_1\|v\|^2_{H^1(\Gamma)}, 
\end{align*}
and \eqref{group} follows easily. 
\end{proof}

 The  proposition below  states  the local well-posedness  of   \eqref{NLS}.
 \begin{proposition}\label{PP1}
For any $u_{0}\in H^{1}(\Gamma)$, there exist $T=T(u_{0})>0$ and a unique solution $u(t)\in C([0,T], H^{1}(\Gamma))\cap C^1([0,T], (H^{1}(\Gamma))')$ of problem \eqref{NLS}.
For each $T_0\in (0, T)$ the mapping
$u_0\in H^1(\Gamma)\mapsto u(t)\in C([0, T_0], H^1(\Gamma))$ is continuous.  Moreover, problem \eqref{NLS} has a maximal solution defined on
an interval of the form  $[0, T_{H^1})$, and the following ``blow-up alternative'' holds: either $T_{H^1} = \infty$ or $T_{H^1} < \infty$ and
$$\lim\limits_{t\to T_{H^1}}\|u(t)\|_{H^1(\Gamma)} =\infty.$$

Finally, the conservation of energy and charge holds:  for $t\in[0, T_{H^1})$
\begin{align}\label{conserv}
E(u(t))=\frac{1}{2} {F}_{\gamma, V}(u(t))-\frac{1}{p+1}\|u(t)\|^{p+1}_{p+1}=E(u_{0}),\quad\left\|u(t)\right\|^2_2=\left\|u_{0}\right\|^2_2.
\end{align}
\end{proposition}
\begin{proof}
A sketch of the  proof was given in \cite{CAFiNo2017}. However,  the rigorous proof (which serves for $p>1$)  might be obtained  repeating the one of \cite[Theorem 4.10.1]{Caz03}. In particular, one needs to  use the fact that $g(u)=|u|^{p-1}u\in C^1(\mathbb{C},\mathbb{C})$ (i.e. $\im(g)$ and $\re(g)$  are $C^1$-functions of  $\re u, \im u$) for $p>1$ and apply inequality \eqref{group}. 

The  proof of   conservation laws \eqref{conserv}  might be obtained involving  regularization procedure analogous to the one introduced in the proof of \cite[Theorem 3.3.5]{Caz03} and using the uniqueness of the solution (see \cite[Theorem 3.3.9]{Caz03}). 
\end{proof}
\begin{remark}
$(i)$\, For $p\geq 4$, the  conservation laws follow easily from Proposition \ref{well_D_H} below and continuous dependence on initial data.  

$(ii)$\,For $1<p<5$, problem  \eqref{NLS} is globally well-posed in $H^1(\Gamma).$ To see that one might  repeat the proof of \cite[Theorem 3.4.1]{Caz03}, where condition (3.4.1)  follows from 
\begin{align*}
    &\|u\|_{p+1}^{p+1}-(Vu,u)_2+\gamma|u_1(0)|^2\leq C\|u'\|_2^{\frac{p-1}{2}}\|u\|_2^{\frac{p+3}{2}}+2\varepsilon\|u'\|_2^2+C_{1}\|u\|_2^2\\&\leq 3\varepsilon\|u'\|_2^2+C_{2}\|u\|_2^{\frac{2(p+3)}{5-p}}+C_{1}\|u\|_2^2\leq 3\varepsilon\|u\|_{H^1(\Gamma)}^2+C(\|u_0\|_2).
\end{align*}
The above estimate is induced by the conservation of charge, estimate \eqref{l_4.9.5}, the Gagliardo-Nirenberg  inequality (see (2.1) in \cite{CAFiNo2017}), and the Young inequality $ab\leq \delta a^q+C_\delta b^{q'}, \,\, \frac{1}{q}+\frac{1}{q'}=1,\,\, q,q'>1, \,\, a,b\geq 0.$ Observe that the key point is that  $q=\frac{4}{p-1}>1$ for $1<p<5.$
\end{remark}
Now, let $m\geq 1+2\omega_0.$ Introduce the norm $\|v\|_{ H_{\gamma, V}}:=\|( H_{\gamma, V} + m)v\|_2$ that endows $\dom( H_{\gamma, V})$ with the structure of a Hilbert space. 
We denote $D_{ H_{\gamma, V}}=(\dom( H_{\gamma, V}), \|\cdot\|_{ H_{\gamma, V}})$.

\begin{proposition}\label{well_D_H}  Let $p\geq 4$ and $u_0\in\dom( H_{\gamma, V})$. Then there exists $T>0$  such that problem \eqref{NLS} has  a unique solution  $u(t)\in C([0,T], D_{ H_{\gamma, V}})\cap C^1([0,T], L^2(\Gamma))$.  Moreover, problem \eqref{NLS} has a maximal solution defined on
an interval of the form  $[0, T_{ H_{\gamma, V}})$, and the following ``blow-up alternative'' holds: either $T_{ H_{\gamma, V}} = \infty$ or $T_{ H_{\gamma, V}} <\infty$ and
$$\lim\limits_{t\to T_{ H_{\gamma, V}}}\|u(t)\|_{ H_{\gamma, V}} =\infty.$$
\end{proposition}
\begin{proof}
The proof repeats the one of \cite[Theorem 2.3]{NataMasa2020} observing that $\dom( H_{\gamma, V})\subset H^1(\Gamma)=\dom(( H_{\gamma, V}+m)^{1/2})$ and, by $m\geq 1+2\omega_0$,
$$\|u\|_\infty\leq C_1\|u\|_{H^1(\Gamma)}\leq C_2\|( H_{\gamma, V}+m)^{1/2}u\|_2\leq C_2 \|( H_{\gamma, V}+m)u\|_2.$$
\end{proof}
\begin{remark}
Notice that due to estimate \eqref{l_4.9.5}, Propositions 2.2 and 2.4 hold for any $\gamma\in\mathbb{R}$ and $V(x)\in (L^1+L^{\infty})(\Gamma).$
\end{remark}
Set
\begin{align*}\label{P}
P(v)=\left\|v'\right\|^{2}_2-\frac{1}{2}\int\limits_{\Gamma}xV'(x)|v(x)|^2\,dx-\frac{\gamma}{2}\left|v_{1}(0)\right|^{2}-\frac{p-1}{2(p+1)}\left\|v\right\|^{p+1}_{p+1}, \quad v\in H^1(\Gamma).
\end{align*}
\begin{proposition}\label{prop-virial} Let $\Sigma(\Gamma)=\{v\in H^1(\Gamma):\, xv\in L^2(\Gamma)\}$. 
Assume that  $u_0\in \Sigma(\Gamma)$, and 
$u(t)$ is the corresponding maximal solution to \eqref{NLS}. 
Then $u(t)\in C([0, T_{H^1}), \Sigma(\Gamma))$, and
the function 
\begin{equation*}\label{well_17}
f(t):=\int\limits_\Gamma x^2|u(t,x)|^2dx=\|xu(t)\|_2^2
\end{equation*}
belongs to $C^2[0,T_{H^1})$. Moreover,
\begin{equation}\label{well_17a}
f'(t)=4\im \int\limits_\Gamma x\overline{u}\partial_x u\,dx,\qquad\text{and }
\end{equation}
\begin{equation}\label{well_18}
f''(t)=8 P(u(t)),\quad t\in [0, T_{H^1}). \quad\text{(virial identity)}
\end{equation}
\end{proposition}
  \begin{proof}
 The proof is similar to the one of \cite[Proposition 6.5.1]{Caz03}. 
 We provide the details since the virial identity  is the key ingredient in the instability analysis. Firstly we show \eqref{well_17a}, secondly we prove \eqref{well_18}
for $u_0\in \dom(H_{\gamma,V})$, then we conclude that \eqref{well_18} holds for $u_0\in H^1(\Gamma)$ using continuous dependence on the initial data.

\textit{Step 1.} \,Let $\varepsilon>0$, define $f_\varepsilon(t)=\|e^{-\varepsilon x^2}xu(t)\|^2_2$, for $t\in [0,T], \,\, T\in (0, T_{H^1})$. Then, observing  that $e^{-2\varepsilon x^2}x^2u(t)\in H^1(\Gamma)$ and  taking $(H^{1})'-H^{1}$ duality product of equation \eqref{NLS} with $ie^{-2\varepsilon x^2}x^2u(t)$,  we get
\begin{equation}\label{well_19}
\begin{split}
f'_\varepsilon(t)&=2\im\int\limits_\Gamma  \left(\partial_xu\,\partial_x(e^{-2\varepsilon x^2}x^2\overline{u})-e^{-2\varepsilon x^2}x^2|u|^{p+1}\right)dx\\&= 4\im \int\limits_\Gamma \left\{e^{-\varepsilon x^2}(1-2\varepsilon x^2)\right\}\overline{u} xe^{-\varepsilon x^2}\partial_xu\,dx.
\end{split}
\end{equation}
Remark that $|e^{-\varepsilon x^2}(1-2\varepsilon x^2)|\leq 2$ for any $x$. From \eqref{well_19}, by the Cauchy-Schwarz inequality, we obtain
\begin{equation}\label{well_20}
\begin{split}
&|f'_\varepsilon(t)|\leq 4\left|\int\limits_\Gamma\left\{e^{-\varepsilon x^2}(1-2\varepsilon x^2)\right\}\overline{u} xe^{-\varepsilon x^2}\partial_xu\, dx\right|\leq 8\int\limits_\Gamma|e^{-\varepsilon x^2}xu\partial_xu|\,dx \\&\leq 8\sum\limits_{j=1}^N\|\partial_xu_j\|_2\|e^{-\varepsilon x^2}xu_j\|_2\leq C\|u\|_{H^1(\Gamma)}\sqrt{f_\varepsilon(t)}.\end{split}
\end{equation}
From \eqref{well_20} one implies
$$\int\limits_0^t\frac{f'_\varepsilon(s)}{\sqrt{f_\varepsilon(s)}}ds\leq C\int\limits_0^t\|u(s)\|_{H^1(\Gamma)}ds,$$
and therefore 
$$\sqrt{f_\varepsilon(t)}\leq \|xu_0\|_2+\frac{C}{2}\int\limits_0^t\|u(s)\|_{H^1(\Gamma)}ds,\,\,t\in [0,T].$$
 Letting $\varepsilon\downarrow 0$ and applying  Fatou's lemma, we get that $xu(t)\in L^2(\Gamma)$ and $f(t)$  is bounded in $[0,T].$
 Observe that from \eqref{well_19} one induces
 \begin{equation}\label{well_21}
  f_\varepsilon(t)=f_\varepsilon(0)+4\im\int\limits_0^t\int\limits_\Gamma \left\{e^{-\varepsilon x^2}(1-2\varepsilon x^2)\right\}\overline{u} xe^{-\varepsilon x^2}\partial_xu\,dx\, ds. \end{equation}
  We have the following estimates for any positive $x$ and $\varepsilon$:
  \begin{equation}\label{well_22}
  \begin{split}
  &e^{-2\varepsilon x^2}x^2|u(t)|^2\leq x^2|u(t)|^2,\\
  &e^{-2\varepsilon x^2}x^2|u_0|^2\leq x^2|u_0|^2,\\
  & |e^{-\varepsilon x^2}(1-2\varepsilon x^2)\overline{u} xe^{-\varepsilon x^2}\partial_xu|\leq 2|\partial_xu\|xu|.
  \end{split}
  \end{equation}
  Having pointwise convergence, and using  \eqref{well_22},  by the Dominated Convergence Theorem we get from \eqref{well_21}
  $$f(t)= \|xu(t)\|_2^2= \|xu_0\|_2^2+4 \im\int\limits_0^t\int\limits_\Gamma x\overline{u}\partial_xu\,dx\,ds.$$
  Since $u(t)$ is strong $H^1$-solution, $f(t)$ is $C^1$-function, and \eqref{well_17a} holds for any $t\in [0, T_{H^1}).$  
  
  Using continuity of $\|xu(t)\|_2$ and the inclusion $u(t)\in C([0, T_{H^1}), H^1(\Gamma))$, by the Brezis-Lieb lemma \cite{LBL}, we get for $t_0, t_n\in [0, T_{H^1})$ $$\lim\limits_{t_n\to t_0}\|xu(t_n)-xu(t_0)\|_2^2=\lim\limits_{t_n\to t_0}\|xu(t_n)\|_2^2-\|xu(t_0)\|_2^2=0,$$ hence $u(t)\in C([0, T_{H^1}), \Sigma(\Gamma)).$
  
  \textit{Step 2. } Let $u_0\in \dom( H_{\gamma, V})$. By Proposition \ref{well_D_H}, the solution $u(t)$ to the corresponding Cauchy problem belongs to $C([0,T_{ H_{\gamma, V}}), D_{ H_{\gamma, V}})\cap C^1([0,T_{ H_{\gamma, V}}), L^2(\Gamma))$.
  
  Let $\varepsilon>0$ and $\theta_\varepsilon(x)=e^{-\varepsilon x^2}$. Define 
    \begin{equation}\label{well_22a}
  h_\varepsilon(t)=\im \int\limits_\Gamma\theta_\varepsilon x\overline{u}\partial_xu\,dx\,\,\,\, \text{for}\,\,\,t\in[0,T],\,\, T\in (0, T_{H_{\gamma, V}}).\end{equation}
  First, let us show that 
  \begin{equation}\label{well_23}
  h'_\varepsilon(t)=-\im\int\limits_\Gamma\partial_tu\left\{2\theta_\varepsilon x\overline{\partial_xu}+(\theta_\varepsilon+x \theta'_\varepsilon)\overline{u}\right\}dx\end{equation} or equivalently
  \begin{equation}\label{well_24}
   h_\varepsilon(t)= h_\varepsilon(0)-\im\int\limits_0^t\int\limits_\Gamma\partial_su\left\{2\theta_\varepsilon x\overline{\partial_xu}+(\theta_\varepsilon+x \theta'_\varepsilon)\overline{u}\right\}dx\,ds.
  \end{equation} 
  Let us prove that identity \eqref{well_24} holds for $u(t)\in C([0,T], H^1(\Gamma))\cap C^1([0,T], L^2(\Gamma)).$ Note that by density argument it is sufficient to show \eqref{well_24} for $u(t)\in C^1([0,T], H^1(\Gamma))\cap C^1([0,T], L^2(\Gamma)).$
  From \eqref{well_22a}, it follows
  \begin{equation}\label{well_24a} h'_\varepsilon(t)=-\im\int\limits_\Gamma\left\{\theta_\varepsilon x\partial_tu\overline{\partial_xu}+\theta_\varepsilon xu\overline{\partial_{xt}^2u}\right\}dx. \end{equation}
  Note that 
  $$\theta_\varepsilon xu\overline{\partial_{xt}^2u}=\theta_\varepsilon xu\overline{\partial_{tx}^2u}=\partial_x\left(\theta_\varepsilon xu\overline{\partial_{t}u}\right)-\theta_\varepsilon u\overline{\partial_{t}u}-\theta_\varepsilon x\partial_xu\overline{\partial_{t}u}-x\theta'_{\varepsilon}u\overline{\partial_{t}u}, $$
  which induces
  $$
\int\limits_{\Gamma}  \theta_\varepsilon xu\overline{\partial_{xt}^2u}\, dx=-\int\limits_{\Gamma} \overline{\partial_{t}u}\left\{\theta_\varepsilon(u+x\partial_xu)+x\theta'_\varepsilon u\right\}dx.
$$
Therefore, from \eqref{well_24a} we get 
$$  h'_\varepsilon(t)=-\im\int\limits_\Gamma\left\{\theta_\varepsilon x\partial_tu\overline{\partial_xu}+\partial_tu\left(\theta_\varepsilon(\overline{u}+x\overline{\partial_xu})+x\theta'_\varepsilon\overline{u}\right)\right\}dx.$$ 
Consequently we obtain \eqref{well_24} for  $u(t)\in C^1([0,T], H^1(\Gamma))\cap C^1([0,T], L^2(\Gamma))$ and hence for $u(t)\in C([0,T], H^1(\Gamma))\cap C^1([0,T], L^2(\Gamma))$ which implies \eqref{well_23}.

 Since $u(t)\in C([0,T_{ H_{\gamma, V}}), D_{H_{\gamma, V}})$, from \eqref{well_23} we get
 \begin{equation}\label{well_25}
  h'_\varepsilon(t)=\re\int\limits_\Gamma( H_{\gamma, V}u-|u|^{p-1}u)\left\{2\theta_\varepsilon x\overline{\partial_xu}+(x\theta_\varepsilon)'\overline{u}\right\}dx.
 \end{equation}
 Below we will consider separately  linear and nonlinear part of identity \eqref{well_25}. 
 Integrating by parts, we obtain 
 \begin{equation}\label{well_26a}
    \begin{split} 
&-\re\int\limits_\Gamma \Delta_\gamma u \left\{2\theta_\varepsilon x\overline{\partial_xu}+(x\theta_\varepsilon)'\overline{u}\right\}dx \\&=-\gamma|u_1(0)|^2+2\int\limits_\Gamma x\theta'_\varepsilon|\partial_xu|^2dx+\int\limits_\Gamma(2\theta'_\varepsilon+x\theta''_\varepsilon)\re(\overline{u}\partial_xu)dx+2\int\limits_\Gamma \theta_\varepsilon|\partial_xu|^2dx.
 \end{split}
 \end{equation}
 Noting that 
 $$\re\left(V(x)u\left\{2\theta_\varepsilon x\overline{\partial_xu}+(x\theta_\varepsilon)'\overline{u}\right\}\right)=\partial_x\left(xV(x)\theta_\varepsilon |u|^2\right)-xV'(x)\theta_\varepsilon |u|^2,$$
 we get 
 \begin{equation}\label{well_26b}
 \re\int\limits_\Gamma V(x)u\left\{2\theta_\varepsilon x\overline{\partial_xu}+(x\theta_\varepsilon)'\overline{u}\right\}dx=-\int\limits_\Gamma xV'(x)\theta_\varepsilon |u|^2 dx.
 \end{equation}
Moreover,
  \begin{equation}\label{well_27} 
  \begin{split}
  & \re\int\limits_\Gamma-|u|^{p-1}u\left\{2\theta_\varepsilon x\overline{\partial_xu}+(x\theta_\varepsilon)'\overline{u}\right\}dx\\&=-\int\limits_\Gamma|u|^{p+1}\theta_\varepsilon\, dx-\int\limits_\Gamma|u|^{p+1}x\theta'_\varepsilon\, dx-\int\limits_\Gamma(|u|^2)^{\frac{p-1}{2}}\partial_x(|u|^2)x\theta_\varepsilon\, dx\\&=-\frac{p-1}{p+1}\int\limits_\Gamma|u|^{p+1}\theta_\varepsilon\, dx-\frac{p-1}{p+1}\int\limits_\Gamma|u|^{p+1}x\theta'_\varepsilon\, dx. 
  \end{split}
  \end{equation}
  Finally, from \eqref{well_25}-\eqref{well_27} we get 
  \begin{equation*}
  \begin{split}
 & h'_\varepsilon(t)=\left[2\int\limits_\Gamma\theta_\varepsilon|\partial_xu|^2dx -\int\limits_\Gamma xV'(x)\theta_\varepsilon |u|^2 dx-\gamma|u_1(0)|^2-\frac{p-1}{p+1}\int\limits_\Gamma|u|^{p+1}\theta_\varepsilon dx\right]\\&+\left[2\int\limits_\Gamma x\theta'_\varepsilon|\partial_xu|^2dx+\int\limits_\Gamma(2\theta'_\varepsilon+x\theta''_\varepsilon)\re(\overline{u}\partial_xu)\,dx\right]-\frac{p-1}{p+1}\int\limits_\Gamma|u|^{p+1}x\theta'_\varepsilon\, dx. 
  \end{split}
  \end{equation*}
  Since $\theta_\varepsilon,\, \theta'_\varepsilon,\, x\theta'_\varepsilon,\, x\theta''_\varepsilon$ are bounded with respect to $x$ and $\varepsilon$, and 
   $$ \theta_\varepsilon\to 1,\,\, \theta'_\varepsilon\to 0,\,\,x\theta'_\varepsilon\to 0,\,\,x\theta''_\varepsilon\to 0\,\, \text{pointwise as}\,\, \varepsilon\downarrow 0, $$
   by the Dominated Convergence Theorem we have
   $$\lim\limits_{\varepsilon\downarrow 0}h'_\varepsilon(t)=2\|\partial_x u\|_2^2-\int\limits_\Gamma xV'(x)|u|^2dx-\gamma|u_1(0)|^2-\frac{p-1}{p+1}\|u\|_{p+1}^{p+1}=:g(t).$$
      Moreover,  again by the Dominated Convergence Theorem, 
      $$\lim\limits_{\varepsilon\downarrow 0}h_\varepsilon(t)=\im \int\limits_\Gamma x\overline{u}\partial_x  u\, dx=:h(t).$$
      Using continuity of $g(t)$ and the fact that the operator $A=\dfrac{d}{dt}$ in  the space $C[0,T]$ with $\dom(A)=C^1[0,T]$  is closed, we arrive at $h'(t)=g(t),\,\, t\in[0,T]$, i.e.
$$ h'(t)= 2\|\partial_x u\|_2^2-\int\limits_\Gamma xV'(x)|u|^2dx-\gamma|u_1(0)|^2-\frac{p-1}{p+1}\|u\|_{p+1}^{p+1},$$
and $h(t)$ is $C^1$ function. Finally, \eqref{well_18} holds for $u_0\in \dom(H_{\gamma, V})$.

\textit{ Step 3.}
 To conclude the proof   consider $\{u_0^n\}_{n\in \mathbb{N}}\subset \dom( H_{\gamma, V})$ such that $u_0^n\to u_0$ in $H^1(\Gamma)$ and $xu_0^n\to xu_0$ in $L^2(\Gamma)$ as $n\to \infty$. Let $u^n(t)$ be the maximal solutions of the corresponding Cauchy problem associated with \eqref{NLS}.
 From \eqref{well_17a} and \eqref{well_18} we obtain 
 $$\|xu^n(t)\|_2^2=\|xu^n_0\|_2^2+4t\im \int\limits_\Gamma x\overline{u_0^n}\partial_xu_0^n\,dx+\int\limits_0^t\int\limits_0^s8P(u^n(y))dy\,ds.$$
Using continuous dependence and repeating the arguments from \cite[Corollary 6.5.3]{Caz03}, we obtain as $n\to \infty$
$$\|xu(t)\|_2^2=\|xu_0\|_2^2+4t\im \int\limits_\Gamma x\overline{u_0}\partial_xu_0\,dx+\int\limits_0^t\int\limits_0^s8P(u(y))dy\,ds,$$ that is \eqref{well_18} holds for  $u_0\in H^1(\Gamma)$.
 \end{proof}

\section{existence of ground states}
\label{S:2}
In this section, we prove  Proposition \ref{GS}. We begin with two technical  lemmas. Throughout this section we assume that $\omega>\omega_0$.
\begin{lemma}\label{V1} If $I_{\omega}(v)<0$, then
\begin{align*}
d_{\omega}<\frac{p-1}{2(p+1)}\left\|v\right\|^{p+1}_{p+1},\,\,\,\mbox{and}\,\,\,d_{\omega}<\frac{p-1}{2(p+1)}\left( {F}_{\gamma, V}(v)+\omega\left\|v\right\|^2_2\right).
\end{align*}
Moreover,
\begin{equation}\label{uu}
\begin{split}
d_{\omega}&=\inf\left\{\frac{p-1}{2(p+1)}\left\|v\right\|^{p+1}_{p+1}\,:\,\,v\in H^{1}(\Gamma)\backslash \{0\},\,\,I_{\omega}(v)\leq0\right\}\\
&=\inf\left\{\frac{p-1}{2(p+1)}\left( {F}_{\gamma, V}(v)+\omega\left\|v\right\|^2_2\right)\,:\,\,v\in H^{1}(\Gamma)\backslash \{0\},\,\,I_{\omega}(v)\leq0\right\}.
\end{split}
\end{equation}
\end{lemma}
\begin{proof}
Noting that 
\begin{align}\label{Somega}
S_{\omega}(v)=\frac{1}{2}I_{\omega}(v)+\frac{p-1}{2(p+1)}\left\|v\right\|^{p+1}_{p+1},\,\,\,\,\,v\in H^{1}(\Gamma),
\end{align}
we get
\begin{align*}
d_{\omega}&=\inf\left\{\frac{p-1}{2(p+1)}\left\|v\right\|^{p+1}_{p+1}\,:\,\,v\in H^{1}(\Gamma)\backslash \{0\},\,\,I_{\omega}(v)=0\right\}.
\end{align*}
Set
\begin{align*}
d^{*}_{\omega}:=\inf\left\{\frac{p-1}{2(p+1)}\left\|v\right\|^{p+1}_{p+1}\,:\,\,v\in H^{1}(\Gamma)\backslash \{0\},\,\,I_{\omega}(v)\leq0\right\}.
\end{align*}
It is clear that $d^{*}_{\omega}\leq d_{\omega}$. Let $v\in H^{1}(\Gamma)\backslash \{0\}$ and $I_{\omega}(v)<0$. Put
\begin{align*}
\lambda_{1}:=\left(\frac{ {F}_{\gamma, V}(v)+\omega\left\|v\right\|^2_2}{\left\|v\right\|^{p+1}_{{p+1}}}\right)^{\frac{1}{p-1}}.
\end{align*}
Then, since $I_{\omega}(\lambda v)=\lambda^{2}\left( {F}_{\gamma, V}(v)+\omega\left\|v\right\|^2_2\right)-\lambda^{p+1}\left\|v\right\|^{p+1}_{p+1}=:f(\lambda)$, we obtain $I_{\omega}(\lambda_{1}v)=0$ and $0<\lambda_{1}<1$ (one needs to remark that $f(1)<0, f(0)=0$, and $f'(\lambda)>0$ for small positive $\lambda$). Hence we have
\begin{align*}
d_{\omega}\leq\frac{p-1}{2(p+1)}\left\|\lambda_{1}v\right\|^{p+1}_{p+1}=\frac{p-1}{2(p+1)}\lambda_{1}^{p+1}\left\|v\right\|^{p+1}_{p+1}<\frac{p-1}{2(p+1)}\left\|v\right\|^{p+1}_{p+1}.
\end{align*}
Thus, we obtain $d_{\omega}\leq d^{*}_{\omega}$. Similarly we can show $d_{\omega}<\frac{p-1}{2(p+1)}\left( {F}_{\gamma, V}(v)+\omega\left\|v\right\|^2_2\right)$ and the second part of \eqref{uu} since we can rewrite
\begin{align*}
d_{\omega}=\inf\left\{\frac{p-1}{2(p+1)}\left( {F}_{\gamma, V}(v)+\omega\left\|v\right\|^2_2\right)\,:\,\,u\in H^{1}(\Gamma)\backslash \{0\},\,\,I_{\omega}(v)=0\right\}.
\end{align*}
\end{proof}
To get the existence of the minimizers of $d_\omega$, one has at a certain point to compare the action $S_{\omega}$ for $\gamma>0$ with the action $S^{0}_{\omega}$ of the nonpotential case ($V(x)\equiv 0$, $\gamma>0$). Set
\begin{equation*}\label{d^0_om}
\begin{split}
S^{0}_{\omega}(v)&=\frac{1}{2}\|v'\|^{2}_2+\frac{\omega}{2}\|v\|_2^2-\frac{\gamma}{2}\left|v_{1}(0)\right|^{2}-\frac{1}{p+1}\|v\|^{p+1}_{p+1},\\
I^{0}_{\omega}(v)&=\|v'\|^{2}_2+\omega\|v\|_2^2-\gamma\left|v_{1}(0)\right|^{2}-\|v\|^{p+1}_{p+1},\\
d^{0}_{\omega}&=\inf\left\{S^{0}_{\omega}(v)\,:\,\,v\in H^{1}(\Gamma)\backslash \{0\},\,\,\,I^{0}_{\omega}(v)=0\right\}\\
&=\inf\left\{\frac{p-1}{2(p+1)}\left\|v\right\|^{p+1}_{p+1}\,:\,\,v\in H^{1}(\Gamma)\backslash \{0\},\,\,\,I^{0}_{\omega}(v)=0\right\},
\end{split}
\end{equation*}
and \begin{align*}\label{M^0_om}
\mathcal{M}^{0}_{\omega}:=\left\{\phi\in H^{1}(\Gamma)\backslash \{0\}\,:\,\,I^{0}_{\omega}(\phi)=0,\,\,S^{0}_{\omega}(\phi)=d^{0}_{\omega}\right\}.
\end{align*}
It is known that  for $\gamma> \gamma^{\ast}$, where $\gamma^*$ is defined by \eqref{gamma_star},  the set $\mathcal{M}^{0}_{\omega}$ is not empty (see \cite{AQFF}). Throughout this section we assume  $\gamma>\gamma^{\ast}$.
\begin{lemma}\label{V2}
$d^{0}_{\omega}>d_{\omega}>0$.
\end{lemma}
\begin{proof} First, we show that $d_{\omega}>0$. Let $v\in H^{1}(\Gamma)\backslash \{0\}$ satisfy $I_{\omega}(v)=0$. Then
\begin{equation*}\label{mm}
\left\|v\right\|^{p+1}_{p+1}= {F}_{\gamma, V}(v)+\omega\left\|v\right\|^2_2.
\end{equation*}
Since $\omega>\omega_0$, by the Sobolev embedding and Lemma \ref{lemma_4.11}-$(ii)$, we have
\begin{align*}
\left\|v\right\|^{2}_{p+1}\leq C_{1}\left\|v\right\|^{2}_{H^{1}(\Gamma)}\leq C_2\left( {F}_{\gamma, V}(v)+\omega\left\|v\right\|^2_2\right)=C_2\left\|v\right\|^{p+1}_{p+1}.
\end{align*}
 Hence we obtain $C_2^{\frac{-1}{p-1}}\leq\left\|v\right\|_{p+1}$. Taking the infimum over $v$, we get $d_{\omega}>0$. Next, we prove $d^{0}_{\omega}>d_{\omega}$. Since $\mathcal{M}^{0}_{\omega}$ is not empty, we can take $\phi\in\mathcal{M}^{0}_{\omega}$. By \textit{Assumption 3}, 
\begin{align*}
I_{\omega}(\phi)=(V\phi,\phi)_2<0.
\end{align*}
Then, by Lemma \ref{V1}, we obtain
\begin{align*}
d_{\omega}<\frac{p-1}{2(p+1)}\left\|\phi\right\|^{p+1}_{p+1}=d^{0}_{\omega}.
\end{align*}
\end{proof}
\begin{lemma}\label{V3}  Let   $\{v_{n}\}\subset H^{1}(\Gamma)\backslash \{0\}$ be a minimizing sequence for $d_{\omega}$, i.e.
$I_{\omega}(v_{n})=0$ and $\lim\limits_{n\rightarrow\infty}S_{\omega}(v_{n})=d_{\omega}$. Then there exist a subsequence $\{v_{n_{k}}\}$ of $\{v_{n}\}$ and $v_{0}\in H^{1}(\Gamma)\backslash \{0\}$ such that $\lim\limits_{k\rightarrow\infty}\left\|v_{n_{k}}-v_{0}\right\|_{H^{1}(\Gamma)}=0$, $I_{\omega}(v_{0})=0$ and $S_{\omega}(v_{0})=d_{\omega}$. Therefore, $\mathcal{M}_{\omega}$ is not empty.
\end{lemma}
\begin{proof} Since $\omega>\omega_0$ and 
\begin{align}\label{mm1}
S_\omega(v_n)=\frac{p-1}{2(p+1)}\left( {F}_{\gamma, V}(v_{n})+\omega\left\|v_{n}\right\|^2_2\right)=\frac{p-1}{2(p+1)}\left\|v_{n}\right\|^{p+1}_{p+1}\underset{n\to\infty}{\longrightarrow} d_{\omega},
\end{align}
 the sequence  $\{v_{n}\}$ is bounded in $H^{1}(\Gamma)$ (see Lemma \ref{lemma_4.11}-$(ii)$). Hence there exist a subsequence $\{v_{n_{k}}\}$ of $\{v_{n}\}$ and $v_{0}\in H^{1}(\Gamma)$ such that $\{v_{n_{k}}\}$ converges weakly to $v_{0}$ in $H^{1}(\Gamma)$. We may assume that  $v_{n_{k}}\neq 0$ and define
\begin{align*}
\lambda_{k}=\left(\frac{\left\|v'_{n_{k}}\right\|^2_2+\omega\left\|v_{n_{k}}\right\|^2_2-\gamma\left|v_{n_{k},1}(0)\right|^{2}}{\left\|v_{n_{k}}\right\|^{p+1}_{p+1}}\right)^{\frac{1}{p-1}}.
\end{align*}
Notice that $\lambda_{k}>0$ and $I^{0}_{\omega}(\lambda_{k}v_{n_{k}})=0$. Therefore, by Lemma \ref{V2} and the definition of $d^{0}_{\omega}$, we obtain
\begin{align}\label{nn1}
d_{\omega}<d^{0}_{\omega}\leq\frac{p-1}{2(p+1)}\left\|\lambda_{k}v_{n_{k}}\right\|^{p+1}_{p+1}=\lambda^{p+1}_{k}\frac{p-1}{2(p+1)}\left\|v_{n_{k}}\right\|^{p+1}_{p+1},\,\,\,\mbox{for all}\,\,k\in \mathbb{N}.
\end{align} Furthermore, by  $I_\omega(v_{n_k})=0$, \eqref{mm1} and the weak continuity of  $(Vv,v)_2=\displaystyle\int\limits_{\Gamma}V(x)|v(x)|^2dx$ (see \cite[Theorem 11.4]{LeiLos01}), we get
\begin{align*}
\lim_{k\rightarrow\infty}\lambda_{k}&=\lim_{k\rightarrow\infty}\left(\frac{\left\|v_{n_{k}}\right\|^{p+1}_{p+1}-(Vv_{n_k},v_{n_k})_2}{\left\|v_{n_{k}}\right\|^{p+1}_{p+1}}\right)^{\frac{1}{p-1}}=\left(\frac{d_{\omega}-\frac{p-1}{2(p+1)}(Vv_{0},v_{0})_2}{d_{\omega}}\right)^{\frac{1}{p-1}}.
\end{align*}
 Taking the limit in \eqref{nn1}, we obtain $d_\omega<\lim\limits_{k\to \infty}\lambda_k^{p+1}d_\omega$.  Since $d_{\omega}>0$, we arrive at $\lim\limits_{k\to \infty}\lambda_k>1$, and consequently  $-(Vv_0,v_0)_2>0$. Thus, $v_{0}\neq0$.\\
  By the weak convergence, we obtain 
\begin{align}\label{Nop}
&\lim_{k\rightarrow\infty}\left\{\left( {F}_{\gamma, V}(v_{n_{k}})- {F}_{\gamma, V}(v_{n_{k}}-v_{0})\right)+\omega\left(\left\|v_{n_{k}}\right\|^2_2-\left\|v_{n_{k}}-v_{0}\right\|^2_2\right)\right\}\\\nonumber
&= {F}_{\gamma, V}(v_{0})+\omega\left\|v_{0}\right\|^2_2.
\end{align}

Next, passing to a subsequence of $\{v_{n_{k}}\}$ if necessary, we may assume that $v_{n_{k}}\underset{k\to\infty}{\longrightarrow} v_{0}$ a.e. on $\Gamma$. 
  Therefore, by the Brezis-Leib lemma \cite{LBL},
$$\lim\limits_{k\rightarrow\infty}I_{\omega}(v_{n_{k}})-I_{\omega}(v_{n_{k}}-v_{0})=\lim\limits_{k\rightarrow\infty}-I_{\omega}(v_{n_{k}}-v_{0}) = I_{\omega}(v_{0}).$$

Since $v_{0}\neq0$, then the right-hand  side of \eqref{Nop} is positive. It follows from \eqref{mm1} and \eqref{Nop} that
\begin{align*}
\frac{p-1}{2(p+1)}\lim_{k\rightarrow\infty}\left( {F}_{\gamma, V}(v_{n_{k}}-v_{0})+\omega\left\|v_{n_{k}}-v_{0}\right\|^2_2\right)
\\<\frac{p-1}{2(p+1)}\lim_{k\rightarrow\infty}\left( {F}_{\gamma, V}(v_{n_{k}})+\omega\left\|v_{n_{k}}\right\|^2_2\right)=d_{\omega}.
\end{align*}
Hence, by \eqref{uu}, we have $I_{\omega}(v_{n_{_{k}}}-v_{0})>0$ for $k$ large enough. Thus, since $-I_{\omega}(v_{n_{k}}-v_{0})\underset{k\to\infty}{\longrightarrow} I_{\omega}(v_{0})$, we obtain $I_{\omega}(v_{0})\leq0$. Then, by \eqref{uu} and the weak lower semicontinuity of norms, we see that
\begin{align*}
d_{\omega}\leq\frac{p-1}{2(p+1)}\left( {F}_{\gamma, V}(v_{0})+\omega\left\|v_{0}\right\|^2_2\right)\leq\frac{p-1}{2(p+1)}\lim_{k\rightarrow\infty}\left( {F}_{\gamma, V}(v_{n_{k}})+\omega\left\|v_{n_{k}}\right\|^2_2\right)=d_{\omega}.
\end{align*}
Therefore, from \eqref{Nop}  we get
\begin{align*}
\lim_{k\rightarrow\infty} {F}_{\gamma, V}(v_{n_{k}}-v_{0})+\omega\left\|v_{n_{k}}-v_{0}\right\|^2_2=0, 
\end{align*}
and consequently, by Lemma \ref{lemma_4.11}-$(ii)$, we have $v_{n_{k}}\underset{k\to\infty}{\longrightarrow} v_{0}$ in $H^{1}(\Gamma)$ and $I_\omega(v_0)=0$. This concludes the proof.
\end{proof}
\begin{proof}[\bf {Proof of Proposition \ref{GS}}]
\textit{Step 1.} We  prove that $\mathcal{G}_{\omega}=\mathcal{M}_{\omega}$.  Let $\varphi\in\mathcal{M}_{\omega}$. Since $I_{\omega}(\varphi)=0$, we have
\begin{align}\label{lli}
\left\langle I^\prime_{\omega}(\varphi),\varphi\right\rangle=2\left( {F}_{\gamma, V}(\varphi)+\omega\left\|\varphi\right\|^2_2\right)-(p+1)\left\|\varphi\right\|^{p+1}_{p+1}=-(p-1)\left\|\varphi\right\|^{p+1}_{p+1}<0.
\end{align} 
There exists a Lagrange multiplier $\mu\in\mathbb{R}$ such that $S^\prime_{\omega}(\varphi)=\mu I^\prime_{\omega}(\varphi)$. Furthermore, since
\begin{align*}
\mu\left\langle I^\prime_{\omega}(\varphi),\varphi\right\rangle=\left\langle S^\prime_{\omega}(\varphi),\varphi\right\rangle= I_{\omega}(\varphi)=0,
\end{align*}
then, by \eqref{lli}, $\mu=0$. Hence $S^\prime_{\omega}(\varphi)=0$. Moreover, for $v\in H^{1}(\Gamma)\backslash \{0\}$ satisfying $S^\prime_{\omega}(v)=0$, we have $I_{\omega}(v)=\left\langle S^\prime_{\omega}(v),v\right\rangle=0$. Then, from the definition of $\mathcal{M}_{\omega}$, we get $S_{\omega}(\varphi)\leq S_{\omega}(v)$. Hence, we obtain $\varphi\in \mathcal{G}_{\omega}$. Now, let $\phi\in \mathcal{G}_{\omega}$. Since $\mathcal{M}_{\omega}$ is not empty, we take $\varphi\in\mathcal{M}_{\omega}$. By the first part of the proof, we have $\varphi\in\mathcal{G}_{\omega}$, therefore $S_{\omega}(\phi)= S_{\omega}(\varphi)=d_{\omega}$. This implies $\phi\in\mathcal{M}_{\omega}$.

\textit{Step 2.} Let $\varphi\in \mathcal{G}_{\omega}$. Below we  show that $\varphi$ has the form $\varphi(x)=e^{i\theta}\phi(x)$ with positive $\phi(x)\in\dom(H_{\gamma, V}).$
 Set $\phi:=\left|\varphi\right|$, then $\left\|\phi^\prime\right\|^2_2\leq\left\|\varphi^\prime\right\|^2_2$ and $S_{\omega}(\phi)\leq S_{\omega}(\varphi)=d_{\omega}$. Using $\mathcal{G}_\omega=\mathcal{M}_\omega$, we obtain $I_{\omega}(\varphi)=0$, then $I_{\omega}(\phi)\leq0$. It follows from Lemma \ref{V1} that $\phi\in\mathcal{M}_{\omega}$ and $S_{\omega}(\varphi)=S_{\omega}(\phi)$. Observe that this implies
\begin{align}\label{lemma1}
\|\phi'\|_2^2=\sum^{N}_{e=1}\int\limits^{\infty}_{0}\left|\phi^\prime_{e}(x)\right|^{2}dx=\sum^{N}_{e=1}\int\limits^{\infty}_{0}\left|\varphi^\prime_{e}(x)\right|^{2}dx=\|\varphi'\|_2^2.
\end{align}

From $S'_\omega(\phi)=0$, repeating the proof of \cite[Theorem 4]{AQFF} (see also \cite[Lemma 4.1]{Ardila2018}),  one gets $\phi\in\dom( H_{\gamma, V})$
and $$ H_{\gamma, V}\phi+\omega\phi-\phi^{p}=0,$$ therefore
\begin{align*}
-\phi^{\prime\prime}_{e}+\omega\phi_{e}+V_e(x)\phi_{e}-\phi_{e}^{p}=0,\quad x\in(0,\infty),\,\, e=1,\ldots, N.
\end{align*}
 Recalling that $V(x)\leq 0$ a.e. on $\Gamma$ (see Remark \ref{V_neg}) and using\cite[Theorem 1]{LVL}, we have that $\phi_{e}$ is either trivial or strictly positive on $(0,\infty)$. Indeed, to prove that,  we need to set  $\beta(s):=\omega s-s^{p}$ and observe that  $\beta(s)\in C^{1}[0,\infty)$ is nondecreasing for $s$ small, and $\beta(0)=\beta(\omega^{\frac{1}{p-1}})=0$. 
 
 Now assume $\phi_{e}(0)=\phi^\prime_{e}(0)=0$ and put
\begin{equation*}
\widetilde{\phi}_{e}(x)=
\begin{cases} 
\phi_{e}(x),\,\,\,\,x\in[0,\infty)\\
0,\,\,\,\,\,x\in(-\delta,0).
\end{cases} 
\end{equation*}
Then, by the Sobolev extension theorem, we have $\widetilde{\phi}_{e}\in H^{2}(-\delta,\infty)$. Moreover,
\begin{equation*}
-\widetilde{\phi}_{e}^{\prime\prime}+\omega\widetilde{\phi}_{e}+V_e(x)\widetilde{\phi}_{e}-\widetilde{\phi}_{e}^{p}=0,\,\,\mbox{on}\,\,(-\delta,\infty).
\end{equation*}
Therefore, by \cite[Theorem 1]{LVL}, arguing as above, we find that $\widetilde{\phi}_{e}=0$ on $(-\delta,\infty)$. 

Next assume $\phi(0)=0$, i.e.  $\phi_1(0)=\ldots=\phi_N(0)$. 
Since $\phi_{e}\in C^{1}(0,\infty)$, $\phi_{e}\geq0$ and $\phi_{e}(0)=0$, then $\phi'_e(0)\geq0$. By  $\sum\limits^{N}_{e=1}\phi'_{e}(0)=-\gamma \phi_{1}(0)=0$, we get $\phi_{e}(0)=\phi^\prime_{e}(0)=0$. Then $\phi_{e}=0$ on $(0,\infty)$ for all $e=1,\ldots,N$, and by continuity $\phi=0$ on $\Gamma$, which is absurd since $\phi\in\mathcal{M}_{\omega}$. Hence $\phi_{e}(0)>0$ for all $e=1,\ldots,N$, therefore $\phi_{e}>0$ on $(0,\infty)$ for all $e=1,\ldots,N$, i.e.  $\phi>0$ on $\Gamma$.

\textit{Step 3.} Now, we can write $\varphi_{e}(x)=\phi_{e}(x)\tau_{e}(x),$ where $\tau_{e}\in C^{1}(0,\infty)$, $\left|\tau_{e}\right|=1$. Then
\begin{equation*}
\varphi^\prime_{e}=\phi^\prime_{e}\tau_{e}+\phi_{e}\tau^\prime_{e}=\tau_{e}(\phi^\prime_{e}+\phi_{e}\overline{\tau}_{e}\tau^\prime_{e}).
\end{equation*}
Using Re$(\overline{\tau}_{e}\tau^\prime_{e})=0$, we have $\left|\varphi^\prime_{e}\right|^{2}=\left|\phi^\prime_{e}\right|^{2}+\left|\phi_{e}\tau^\prime_{e}\right|^{2}$. Therefore, from \eqref{lemma1} we obtain
\begin{equation*}
\sum^{N}_{e=1}\int\limits^{\infty}_{0}\left|\phi^\prime_{e}\right|^{2}dx=\sum^{N}_{e=1}\int\limits^{\infty}_{0}\left|\varphi^\prime_{e}\right|^{2}dx=\sum^{N}_{e=1}\int\limits^{\infty}_{0}\left|\phi^\prime_{e}\right|^{2}dx+\sum^{N}_{e=1}\int\limits^{\infty}_{0}\left|\phi_{e}\tau^\prime_{e}\right|^{2}dx.
\end{equation*}
So far as $\phi_{e}>0$, we have $\tau^\prime_{e}=0$ for all $e=1,\ldots,N$. Since $\tau_{e}\in C^{1}(0,\infty)$, there exists a constant $\theta_{e}\in\mathbb{R}$ such that $\tau_{e}(x)=e^{i\theta_{e}}$ on $(0,\infty)$. By the continuity at the vertex, we obtain $\theta_{e}=\theta=const$  for all $e=1,\ldots,N$. This ends the proof.

 $\mathrm{Re}(\overline{\tau_e}\tau'_e)$
\end{proof}

\section{instability of standing waves}
\label{S:3}
In this section, we prove Theorem \ref{INLS} and Corollary  \ref{CINLS}. 
\subsection{Proof of the  main result}\label{main_res}
We begin with the following lemma.
\begin{lemma}\label{T1}
Let $\varphi_{\omega}\in\mathcal{M}_{\omega}$. Then
\begin{align*}
{(i)}\,\,\, \left\|\varphi_{\omega}\right\|^{p+1}_{p+1}&=\inf\left\{\left\|v\right\|^{p+1}_{p+1}\,:\,\,v\in H^{1}(\Gamma)\backslash \{0\},\,\,I_{\omega}(v)=0\right\}\\
&=\inf\left\{\left\|v\right\|^{p+1}_{p+1}\,:\,\,v\in H^{1}(\Gamma)\backslash \{0\},\,\,I_{\omega}(v)\leq0\right\},\\
{(ii)}\quad S_{\omega}(\varphi_{\omega})&=\inf\{S_{\omega}(v)\,:\,\,v\in H^{1}(\Gamma),\,\,\left\|v\right\|^{p+1}_{p+1}=\left\|\varphi_{\omega}\right\|^{p+1}_{p+1}\}.
\end{align*}
\end{lemma}
\begin{proof}$(i)$  This is an immediate consequence of Lemma \ref{V1}\\
$(ii)$ Set $d^{**}_{\omega}:=\inf\{S_{\omega}(v)\,:\,\,v\in H^{1}(\Gamma),\,\,\left\|v\right\|^{p+1}_{p+1}=\left\|\varphi_{\omega}\right\|^{p+1}_{p+1}\}$. As far as $d^{**}_{\omega}\leq S_{\omega}(\varphi_{\omega})$, it suffices to prove $S_{\omega}(\varphi_{\omega})\leq d^{**}_{\omega}$. If $v\in H^{1}(\Gamma)$ satisfies $\left\|v\right\|^{p+1}_{p+1}=\left\|\varphi_{\omega}\right\|^{p+1}_{p+1}$, then, by item $(i)$ and \eqref{Somega}, we have $I_{\omega}(v)\geq0$. Hence, by \eqref{Somega},
\begin{align*}
S_{\omega}(\varphi_{\omega})=\frac{p-1}{2(p+1)}\left\|\varphi_{\omega}\right\|^{p+1}_{p+1}=\frac{p-1}{2(p+1)}\left\|v\right\|^{p+1}_{p+1}\leq S_{\omega}(v).
\end{align*}
Thus, we obtain $S_{\omega}(\varphi_{\omega})\leq d^{**}_{\omega}$.
\end{proof}
Recall that $$P(v)=\left\|v'\right\|^{2}_2-\frac{1}{2}\int\limits_{\Gamma}xV'(x)|v(x)|^2\,dx-\frac{\gamma}{2}\left|v_{1}(0)\right|^{2}-\frac{p-1}{2(p+1)}\left\|v\right\|^{p+1}_{p+1}.$$
\begin{lemma}\label{T2}
If $\partial^{2}_{\lambda}E(\varphi^{\lambda}_{\omega})\left|_{\lambda=1}\right.<0$, then there  exist $\delta>0$ and $\varepsilon >0$ such that the following holds: for any $v\in N_{\varepsilon}(\varphi_{\omega})$ satisfying $\left\|v\right\|^2_2\leq\left\|\varphi_{\omega}\right\|^2_2$, there exists $\lambda(v)\in(1-\delta,1+\delta)$ such that $E(\varphi_{\omega})\leq E(v)+(\lambda(v)-1)P(v)$, where $N_{\varepsilon}(\varphi_{\omega})$ is  defined by  \eqref{N}.
\end{lemma}
\begin{proof}Since $\partial^{2}_{\lambda}E(\varphi^{\lambda}_{\omega})\left|_{\lambda=1}\right.<0$ and $\partial^{2}_{\lambda}E(v^{\lambda})$ is continuous in $v$  (we mean "orbit"-continuity) and $\lambda$, there exist positive constants $\varepsilon$ and $\delta$ such that $\partial^{2}_{\lambda}E(v^{\lambda})<0$ for any $v\in N_{\varepsilon}(\varphi_{\omega})$ and $\lambda\in(1-\delta,1+\delta)$. Using $P(v)=\partial_{\lambda}E(v^{\lambda})\left|_{\lambda=1}\right.$, the Taylor expansion at $\lambda=1$ gives
\begin{align}\label{Eq5}
E(v^{\lambda})\leq E(v)+(\lambda-1)P(v),\quad\lambda\in(1-\delta,1+\delta),\quad v\in N_{\varepsilon}(\varphi_{\omega}).
\end{align}
Let $v\in N_{\varepsilon}(\varphi_{\omega})$ satisfy $\left\|v\right\|^2_2\leq\left\|\varphi_{\omega}\right\|^2_2$. We define
\begin{align*}
\lambda(v):=\left(\frac{\left\|\varphi_{\omega}\right\|^{p+1}_{p+1}}{\left\|v\right\|^{p+1}_{p+1}}\right)^{\frac{2}{p-1}}.
\end{align*}
Then, $\left\|v^{\lambda(v)}\right\|^{p+1}_{p+1}=\left\|\varphi_{\omega}\right\|^{p+1}_{p+1}$ and we can take $\varepsilon$ small enough to guarantee
$\lambda(v)\in(1-\delta,1+\delta)$. Since $\left\|v^{\lambda(v)}\right\|^2_2=\left\|v\right\|^2_2\leq\left\|\varphi_{\omega}\right\|^2_2$, by Lemma \ref{T1}-$(ii)$, we have
\begin{align*} 
E(v^{\lambda(v)})=S_{\omega}(v^{\lambda(v)})-\frac{\omega}{2}\left\|v^{\lambda(v)}\right\|^2_2\geq S_{\omega}(\varphi_{\omega})-\frac{\omega}{2}\left\|\varphi_{\omega}\right\|^2_2=E(\varphi_{\omega}),
\end{align*}
which together with \eqref{Eq5} implies that $E(\varphi_{\omega})\leq E(v)+(\lambda(v)-1)P(v)$.
\end{proof}
To prove Theorem \ref{INLS}, we introduce the following definition.
\begin{definition}
Let $\varepsilon$ be the positive constant given by Lemma \ref{T2}. Set
\begin{align*}
\mathcal{Z}_{\varepsilon}(\varphi_{\omega}):=\{v\in N_{\varepsilon}(\varphi_{\omega})\,:\,\,E(v)<E(\varphi_{\omega}),\,\,\left\|v\right\|^2_2\leq\left\|\varphi_{\omega}\right\|^2_2,\,\,P(v)<0\},
\end{align*}
and for any $u_{0}\in N_{\varepsilon}(\varphi_{\omega})$, we define the exit time from $N_{\varepsilon}(\varphi_{\omega})$ by
\begin{align*}
T_{\varepsilon}(u_{0})=\sup\{T>0\,:\,\,u(t)\in N_{\varepsilon}(\varphi_{\omega}),\,\,0\leq t\leq T\},
\end{align*}
with $u(t)$ being a solution of \eqref{NLS}. 
\end{definition}
\begin{lemma}\label{T3}
Assume $\partial^{2}_{\lambda}E(\varphi^{\lambda}_{\omega})\left|_{\lambda=1}\right.<0$, then for any $u_{0}\in \mathcal{Z}_{\varepsilon}(\varphi_{\omega})$, there exists $b=b(u_{0})>0$ such that $P(u(t))\leq-b$ for $0\leq t< T_{\varepsilon}(u_{0})$.
\end{lemma}	
\begin{proof}Set $b_0:=E(\varphi_{\omega})-E(u_{0})>0$, with $u_{0}\in\mathcal{Z}_{\varepsilon}(\varphi_{\omega})$. From the conservation of energy  and Lemma \ref{T2}, we have
\begin{align}\label{Eq7}
b_0\leq(\lambda(u(t))-1)P(u(t)),\quad 0\leq t< T_{\varepsilon}(u_{0}).
\end{align}
Therefore, for $0\leq t< T_{\varepsilon}(u_{0})$ we get $P(u(t))\neq0$. Indeed, if $P(u(t_0))=0$ for some $t_0\in[0, T(u_0)),$ then from \eqref{Eq7} it follows $b_0\leq 0$, which contradicts the definition of $b_0$.  Since $P(u_{0})<0$ and the function $t\mapsto P(u(t))$ is continuous, we see that $P(u(t))<0$ for $0\leq t< T_{\varepsilon}(u_{0})$ and hence $\lambda(u(t))-1<0$ for $0\leq t< T_{\varepsilon}(u_{0})$. Thus, from Lemma \ref{T2} and \eqref{Eq7}, we have
\begin{align*}
P(u(t))\leq\frac{b_0}{\lambda(u(t))-1}\leq\frac{-b_0}{\delta},\quad 0\leq t< T_{\varepsilon}(u_{0}).
\end{align*}
Hence, taking $b=\dfrac{b_0}{\delta}$, we arrive at $P(u(t))\leq-b$ for $0\leq t< T_{\varepsilon}(u_{0})$.
\end{proof}
Now we are ready to prove  Theorem \ref{INLS}.
\begin{proof}[\bf {Proof of Theorem \ref{INLS}}]
Observe that $P(v)=\partial_{\lambda}S_{\omega}(v^{\lambda})|_{\lambda=1}= \left\langle S^\prime_{\omega}(v),\partial_{\lambda}v^{\lambda}|_{\lambda=1}\right\rangle$.  Since $S^\prime_{\omega}(\varphi_{\omega})=0$, we obtain $P(\varphi_{\omega})=\partial_{\lambda}S_{\omega}(\varphi^{\lambda}_{\omega})\left|_{\lambda=1}\right.=0$.   Moreover, by $P(\varphi^{\lambda}_{\omega})=\lambda\partial_{\lambda}E(\varphi^{\lambda}_{\omega})$, we have $\partial_{\lambda}E(\varphi^{\lambda}_{\omega})\left|_{\lambda=1}\right.=0$. Then, from the assumption $\partial^{2}_{\lambda}E(\varphi^{\lambda}_{\omega})\left|_{\lambda=1}\right.<0$, we get $E(\varphi^{\lambda}_{\omega})<E(\varphi_{\omega})$ and $P(\varphi^{\lambda}_{\omega})<0$ for $\lambda>1$ close enough to $1$.

Let $\varepsilon>0$ be given by Lemma \ref{T2}.  Since 
$\lim\limits_{\lambda\rightarrow1}\left\|\varphi^{\lambda}_{\omega}-\varphi_{\omega}\right\|_{H^{1}(\Gamma)}=0$ and $\left\|\varphi^{\lambda}_{\omega}\right\|^2_2=\left\|\varphi_{\omega}\right\|^2_2$, by continuity of $E$ and $P$, for any $\delta\leq \varepsilon$ there exists $\lambda_1$ such that $\varphi^{\lambda_1}_{\omega}\in \mathcal{Z}_{\frac{\delta}{2}}(\varphi_\omega).$ 

 Suppose that   $\chi\in C_{c}^{\infty}(\mathbb{R}^{+})$ is the  function satisfying 
 $$0\leq\chi\leq1,\quad  \chi(x)=1,\,\,\text{if}\,\, x\in[0,1],\quad \text{and}\,\, \chi(x)=0\,\, \text{if}\,\, x\geq2.$$ For $a>0$, we define  $\chi_{a}\in C^{\infty}_{c}(\Gamma)$ by
\begin{align*}
(\chi_{a})_{e}(x)=\chi\left(\frac{x}{a}\right),\quad x\in\mathbb{R}^{+},\,\,\,e=1,\ldots,N.
\end{align*}
Then  we have $\lim\limits_{a\rightarrow\infty}\left\|\chi_{a}\varphi^{\lambda_{1}}_{\omega}-\varphi^{\lambda_{1}}_{\omega}\right\|_{H^{1}(\Gamma)}=0$ and $\left\|\chi_{a}\varphi^{\lambda_{1}}_{\omega}\right\|^2_2\leq \left\|\varphi^{\lambda_{1}}_{\omega}\right\|^2_2=\left\|\varphi_{\omega}\right\|^2_2$ for all $a>0$. Thus, by continuity of $E$ and $P$, for any $\delta\leq \varepsilon$ there exists $a_{1}>0$ such that  $\chi_{a_{1}}\varphi^{\lambda_{1}}_{\omega}\in\mathcal{Z}_{\frac{\delta}{2}}(\varphi_\omega^{\lambda_1})$, therefore  $\chi_{a_{1}}\varphi^{\lambda_{1}}_{\omega}\in\mathcal{Z}_{\delta}(\varphi_\omega)\subseteq\mathcal{Z}_{\varepsilon}(\varphi_\omega)$.

Observe that  $\chi_{a_{1}}\varphi^{\lambda_{1}}_{\omega}\in \Sigma(\Gamma)$ (see Proposition \ref{prop-virial} for the definition of $\Sigma(\Gamma)$),
and by  virial identity \eqref{well_18}, we see that
\begin{align}\label{EEE}
\frac{d^{2}}{dt^{2}}\left\|xu_{1}(t)\right\|^2_2=8P(u_{1}(t)),\qquad 0\leq t\leq T_\varepsilon(\chi_{a_{1}}\varphi^{\lambda_{1}}_{\omega}),
\end{align}
where $u_{1}(t)$ is the solution to \eqref{NLS} with $u_{1}(0)=\chi_{a_{1}}\varphi^{\lambda_{1}}_{\omega}$. From Lemma \ref{T3}, there exists $b=b(\lambda_{1},a_{1})>0$ such that
\begin{align}\label{Ex}
P(u_1(t))\leq-b,\qquad 0\leq t< T_\varepsilon(\chi_{a_{1}}\varphi^{\lambda_{1}}_{\omega}).
\end{align}
Then, from \eqref{Ex} and \eqref{EEE}, we can see that $T_\varepsilon(\chi_{a_{1}}\varphi^{\lambda_{1}}_{\omega})<\infty$.

Summarizing the above, we affirm: there exists $\varepsilon>0$ (given by Lemma \ref{T2}) such that for all $\delta>0$ there exists $u_0=\chi_{a_{1}}\varphi^{\lambda_{1}}_{\omega}\in N_\delta(\varphi_\omega)$ and $t_1>0$ such that the corresponding solution $u_1(t)$ of \eqref{NLS} satisfies $u_1(t_1)\notin N_\varepsilon(\varphi_\omega). $ Hence, the standing wave solution $e^{i\omega t}\varphi_{\omega}$ of \eqref{NLS} is orbitally unstable. 
\end{proof}
\subsection{Rescaled variational problem and proof of Corollary \ref{CINLS}}\label{proof_cor}
Assume that $V(x)=\dfrac{-\beta}{x^\alpha},\, \beta>0,\, 0<\alpha<1$. Recall that $v^{\lambda}(x)=\lambda^{1/2}v(\lambda x)$ for $\lambda>0$. By simple computations, we have
\begin{align*}
E(v^{\lambda})=\frac{\lambda^{2}}{2}\|v'\|^{2}_2+\frac{\lambda^{\alpha}}{2}(Vv,v)_2-\frac{\lambda}{2}\gamma\left|v_{1}(0)\right|^{2}-\frac{\lambda^{\frac{p-1}{2}}}{p+1}\|v\|^{p+1}_{p+1},
\end{align*}
\begin{align*}
\partial^{2}_{\lambda}E(v^{\lambda})\left|_{\lambda=1}\right.=\|v'\|^{2}_2+\frac{\alpha(\alpha-1)}{2}(Vv,v)_2-\frac{(p-1)(p-3)}{4(p+1)}\|v\|^{p+1}_{p+1}.
\end{align*}
Since $P(\varphi_{\omega})=\partial_{\lambda}S_{\omega}(\varphi^{\lambda}_{\omega})\left|_{\lambda=1}\right.=0$, then we get 
\begin{align*}
\partial^{2}_{\lambda}E(\varphi^{\lambda}_{\omega})\left|_{\lambda=1}\right.=-\frac{\alpha(2-\alpha)}{2} (V\varphi_\omega,\varphi_\omega)_2+\frac{\gamma}{2}\left|\varphi_{\omega,1}(0)\right|^{2}-\frac{(p-1)(p-5)}{4(p+1)}\|\varphi_{\omega}\|^{p+1}_{p+1},
\end{align*}
and $\partial^{2}_{\lambda}E(\varphi^{\lambda}_{\omega})\left|_{\lambda=1}\right.<0$ is equivalent to
\begin{align}\label{d2E}
\frac{-\alpha(2-\alpha) (V\varphi_\omega,\varphi_\omega)_2+\gamma\left|\varphi_{\omega,1}(0)\right|^{2}}{\|\varphi_{\omega}\|^{p+1}_{p+1}}<\frac{(p-1)(p-5)}{2(p+1)}.
\end{align}
Below we  prove that the left-hand side of \eqref{d2E} converges to $0$ as $\omega\rightarrow\infty$. To this end, we consider the following rescaling of $\varphi_{\omega}\in \mathcal{M}_{\omega}$:
\begin{align}\label{res}
\varphi_{\omega}(x)=\omega^{\frac{1}{p-1}}\widetilde{\varphi}_{\omega}(\sqrt{\omega}x),\quad\omega\in(\omega_0,\infty),
\end{align}
and observe
\begin{equation}\label{E1}
\begin{split}
&\frac{-\omega^{-\frac{2-\alpha}{2}}\alpha(2-\alpha)(V\widetilde{\varphi}_\omega, \widetilde{\varphi}_\omega)_2+\omega^{-\frac{1}{2}}\gamma\left|\widetilde{\varphi}_{\omega,1}(0)\right|^{2}}{\|\widetilde{\varphi}_{\omega}\|^{p+1}_{p+1}}\\
&=\frac{-\alpha(2-\alpha)(V\varphi_\omega,\varphi_\omega)_2+\gamma\left|\varphi_{\omega,1}(0)\right|^{2}}{\|\varphi_{\omega}\|^{p+1}_{p+1}}.
\end{split}
\end{equation}
Put
\begin{align*}
\widetilde{I}_{\omega}(v):&=\|v'\|^{2}_2+\|v\|^{2}_2-\omega^{-\frac{2-\alpha}{2}}\beta\int\limits_{\Gamma}\frac{\left|v(x)\right|^{2}}{x^\alpha}dx-\omega^{-\frac{1}{2}}\gamma\left|v_{1}(0)\right|^{2}-\|v\|^{p+1}_{p+1},\\
\widetilde{I}_0(v):&=\|v'\|^{2}_2+\|v\|^{2}_2-\|v\|^{p+1}_{p+1}.
\end{align*}
Consider the minimization problem
\begin{align}\label{d10}
\widetilde{d}_0:=\inf\left\{\left\|v\right\|^{p+1}_{p+1}\,:\,\,v\in H^{1}(\Gamma)\backslash \{0\},\,\,\widetilde{I}_0(v)\leq0 \right\}.
\end{align}
In \cite[Theorem 3]{AQFF} it was shown that   $\widetilde{d}_0>0$.
The following  lemma is the key result to  prove Corollary \ref{CINLS}.

\begin{lemma}\label{L1}
Assume $\gamma>0$, $\beta>0$, $0<\alpha<1$ and $p>5$. Let $\varphi_{\omega}\in\mathcal{M}_{\omega}$, and  $\widetilde{\varphi}_{\omega}(x)$ be the rescaled function given in \eqref{res}. Then 
\begin{align*}
(i)\, &\lim_{\omega\rightarrow\infty}\left\|\widetilde{\varphi}_{\omega}\right\|^{p+1}_{p+1}=\widetilde{d}_0,\hspace{10cm}\left.\right.\\
(ii)\, &\lim_{\omega\rightarrow\infty}\widetilde{I}_0(\widetilde{\varphi}_{\omega})=0,\hspace{11cm}\left.\right.\\
(iii)\, &\lim_{\omega\rightarrow\infty}\left\|\widetilde{\varphi}_{\omega}\right\|^{2}_{H^{1}(\Gamma)}=\widetilde{d}_0.\hspace{10.2cm}\left.\right.
\end{align*}
\end{lemma}
\begin{proof}
Notice that 
\begin{equation}\label{Min2}
\begin{split}
&\left\|\widetilde{\varphi}_{\omega}\right\|^{p+1}_{p+1}=\inf\left\{\left\|v\right\|^{p+1}_{p+1}\,:\,\,v\in H^{1}(\Gamma)\backslash\{0\},\,\widetilde{I}_{\omega}(v)=0\right\}\\
&=\inf\left\{\left\|v\right\|^{p+1}_{p+1}\,:\,\,v\in H^{1}(\Gamma)\backslash\{0\},\,\widetilde{I}_{\omega}(v)\leq0\right\}:=\widetilde{d}_\omega.
\end{split}
\end{equation}
By definition we have 
\begin{equation}\label{I_1}
\widetilde{I}_0(v)=\widetilde{I}_{\omega}(v)
-\omega^{-\frac{2-\alpha}{2}}(Vv, v)_2+\omega^{-\frac{1}{2}}\gamma\left|v_1(0)\right|^{2},\quad\text{and} 
\end{equation}
\begin{equation}\label{I_2}
\widetilde{I}_0(v)=\lambda^{-2}\widetilde{I}_0(\lambda v)+(\lambda^{p-1}-1)\|v\|_{p+1}^{p+1}.
\end{equation}
Using, \eqref{I_1}, \eqref{I_2}, $\widetilde{I}_{\omega}(\widetilde{\varphi}_{\omega})=0$, estimate  \eqref{ee2}, and the Sobolev embedding, for any $\lambda>1$ we  get
\begin{equation}\label{I_tild}
\begin{split}
\lambda^{-2}\widetilde{I}_0(\lambda\widetilde{\varphi}_{\omega})&=-\omega^{-\frac{2-\alpha}{2}}(V\widetilde{\varphi}_\omega, \widetilde{\varphi}_\omega)_2+\omega^{-\frac{1}{2}}\gamma\left|\widetilde{\varphi}_{\omega,1}(0)\right|^{2}-(\lambda^{p-1}-1)\|\widetilde{\varphi}_{\omega}\|^{p+1}_{p+1}\\
&\leq C_1\omega^{-\frac{2-\alpha}{2}}\left\|\widetilde{\varphi}_{\omega}\right\|^{2}_{H^{1}(\Gamma)}+C_2\omega^{-\frac{1}{2}}\gamma\left\|\widetilde{\varphi}_{\omega}\right\|^{2}_{H^{1}(\Gamma)}-(\lambda^{p-1}-1)\|\widetilde{\varphi}_{\omega}\|^{p+1}_{p+1}.
\end{split}
\end{equation}
Moreover, from $\widetilde{I}_{\omega}(\widetilde{\varphi}_{\omega})=0$, we deduce
\begin{align*}
\left\|\widetilde{\varphi}_{\omega}\right\|^{2}_{H^{1}(\Gamma)}&=-\omega^{-\frac{2-\alpha}{2}}(V\widetilde{\varphi}_\omega, \widetilde{\varphi}_\omega)_2+\omega^{-\frac{1}{2}}\gamma\left|\widetilde{\varphi}_{\omega,1}(0)\right|^{2}+\|\widetilde{\varphi}_{\omega}\|^{p+1}_{p+1}\\
&\leq C_1\omega^{-\frac{2-\alpha}{2}}\left\|\widetilde{\varphi}_{\omega}\right\|^{2}_{H^{1}(\Gamma)}+C_2\omega^{-\frac{1}{2}}\gamma\left\|\widetilde{\varphi}_{\omega}\right\|^{2}_{H^{1}(\Gamma)}+\|\widetilde{\varphi}_{\omega}\|^{p+1}_{p+1}.
\end{align*}
This implies
\begin{align*}
\left(1-C_1\omega^{-\frac{2-\alpha}{2}}-C_2\omega^{-\frac{1}{2}}\gamma\right)\left\|\widetilde{\varphi}_{\omega}\right\|^{2}_{H^{1}(\Gamma)}\leq \|\widetilde{\varphi}_{\omega}\|^{p+1}_{p+1}.
\end{align*}
Since for  $\omega$ sufficiently large $\left(1-C_1\omega^{-\frac{2-\alpha}{2}}-C_2\omega^{-\frac{1}{2}}\gamma\right)>0$, from \eqref{I_tild} we get 
\begin{align}\label{Id0}
\lambda^{-2}\widetilde{I}_0(\lambda\widetilde{\varphi}_{\omega})\leq-\left(\lambda^{p-1}-1-\frac{C_1\omega^{-\frac{2-\alpha}{2}}+C_2\omega^{-\frac{1}{2}}\gamma}{1-C_1\omega^{-\frac{2-\alpha}{2}}-C_2\omega^{-\frac{1}{2}}\gamma}\right)\|\widetilde{\varphi}_{\omega}\|^{p+1}_{p+1}.
\end{align}
Hence for any $\lambda>1$, there exists $\omega_1=\omega_1(\lambda)\in(\omega_0,\infty)$ such that $\widetilde{I}_0(\lambda\widetilde{\varphi}_{\omega})<0$ for  $\omega\in(\omega_1,\infty)$. Thus, by \eqref{d10}, $\widetilde{d_0}\leq \lambda^{p+1}\|\widetilde{\varphi}_\omega\|_{p+1}^{p+1}$ for $\omega\in(\omega_1,\infty)$. Observe that $\widetilde{I}_0(v)\leq 0$ implies $\widetilde{I}_\omega(v)\leq 0$, then from \eqref{Min2} we obtain $\widetilde{d}_\omega=\|\widetilde{\varphi}_\omega\|_{p+1}^{p+1}\leq\widetilde{d_0}$. Therefore,
\begin{equation}\label{d_0-i}
\lambda^{-(p+1)}\widetilde{d}_0\leq \|\widetilde{\varphi}_{\omega}\|^{p+1}_{p+1}\leq\widetilde{d}_0,\quad\omega\in(\omega_1,\infty).
\end{equation}
Letting $\lambda\downarrow 1$, we get that $\omega\to \infty$, and from \eqref{d_0-i} it follows  $(i)$.

Now,  assume that  $\lambda=1$ in  \eqref{Id0}, then using  $(i)$, we deduce
\begin{align}\label{sup}
\limsup_{\omega\rightarrow\infty}\widetilde{I}_0(\widetilde{\varphi}_{\omega})\leq0.
\end{align}
Furthermore, define
$$\lambda_1(\omega)=\left(\frac{\|\widetilde{\varphi}'_\omega\|_2^2+\|\widetilde{\varphi}_\omega\|_2^2}{\|\widetilde{\varphi}_\omega\|_{p+1}^{p+1}}\right)^{\frac 1{p-1}}>0,$$ then $\widetilde{I}_0(\lambda_1(\omega)\widetilde{\varphi}_{\omega})=0$.  Therefore, we have
\begin{align}\label{nj}
\widetilde{d}_0\leq\lambda_1(\omega)^{p+1} \|\widetilde{\varphi}_{\omega}\|^{p+1}_{p+1}.
\end{align}
Thus, by $(i)$ and \eqref{nj}, we arrive at
\begin{align*}
\liminf_{\omega\rightarrow\infty}\lambda_1(\omega)\geq\liminf_{\omega\rightarrow\infty}\left(\frac{\widetilde{d}_0}{\left\|\widetilde{\varphi}_{\omega}\right\|^{p+1}_{p+1}}\right)^{\frac{1}{p+1}}=1.
\end{align*}
Moreover, by \eqref{I_2},
$\widetilde{I}_0(\lambda_1(\omega)\widetilde{\varphi}_{\omega})=0$ and $(i)$, we have
\begin{align*}
\liminf_{\omega\rightarrow\infty}\widetilde{I}_0(\widetilde{\varphi}_{\omega})=\liminf_{\omega\rightarrow\infty}(\lambda_1(\omega)^{p-1}-1)\left\|\widetilde{\varphi}_{\omega}\right\|^{p+1}_{p+1}\geq0,
\end{align*}
which together with \eqref{sup} implies $(ii)$. Finally, from $(i)$ and $(ii)$, we obtain
\begin{align*}
\widetilde{d}_0=\lim_{\omega\rightarrow\infty}\left\|\widetilde{\varphi}_{\omega}\right\|^{p+1}_{p+1}=\lim_{\omega\rightarrow\infty}\left\|\widetilde{\varphi}_{\omega}\right\|^{2}_{H^{1}(\Gamma)},
\end{align*}
which shows $(iii)$. 
\end{proof}

\begin{proof}[\bf {Proof of Corollary \ref{CINLS}}] Recall that, by Theorem \ref{INLS}, if $\partial_{\lambda}^{2} E\left(\varphi_{\omega}^{\lambda}\right)|_{\lambda=1}<0$, then $e^{it}\varphi_{\omega}(x)$ is orbitally unstable.
Since $$\partial_{\lambda}^{2} E\left(\varphi_{\omega}^{\lambda}\right)|_{\lambda=1}<0 \quad \Longleftrightarrow \quad \frac{-\alpha(2-\alpha)\left(V \varphi_{\omega}, \varphi_{\omega}\right)_{2}+\gamma\left|\varphi_{\omega, 1}(0)\right|^{2}}{\left\|\varphi_{\omega}\right\|_{p+1}^{p+1}}<\frac{(p-1)(p-5)}{2(p+1)}, 
$$
 by \eqref{E1}, it suffices to prove
\begin{align}\label{E2}
\lim_{\omega\rightarrow\infty}\frac{-\omega^{-\frac{2-\alpha}{2}}\alpha(2-\alpha)(V\widetilde{\varphi}_{\omega}, \widetilde{\varphi}_{\omega})_2+\omega^{-\frac{1}{2}}\gamma\left|\widetilde{\varphi}_{\omega,1}(0)\right|^{2}}{\|\widetilde{\varphi}_{\omega}\|^{p+1}_{p+1}}=0.
\end{align}
We have
\begin{align*}
0\leq -\omega^{-\frac{2-\alpha}{2}}\alpha(2-\alpha)(V\widetilde{\varphi}_{\omega}, \widetilde{\varphi}_{\omega})_2+\omega^{-\frac{1}{2}}\gamma\left|\widetilde{\varphi}_{\omega,1}(0)\right|^{2}\\
\leq \left(C_1\omega^{-\frac{2-\alpha}{2}}+C_2\omega^{-\frac{1}{2}}\gamma\right)\left\|\widetilde{\varphi}_{\omega}\right\|^{2}_{H^{1}(\Gamma)}.
\end{align*}
Hence, by Lemma \ref{L1}-$(i)$,$(iii)$, we obtain \eqref{E2}. This concludes the proof.
\end{proof}
\subsection{Instability results in $H_{\eq}^1(\Gamma)$}\label{subsec_eq}
We discuss counterparts of Proposition \ref{GS}, Theorem \ref{INLS}, Corollary \ref{CINLS} for  arbitrary $\gamma\in\mathbb{R}$ and symmetric $V(x)$, i.e. $V_1(x)=\ldots=V_N(x)$, in the space 
\begin{equation*}\label{eq}
  H_{\eq}^1(\Gamma)=\{v\in H^1(\Gamma):\,v_1(x)=\ldots=v_N(x),\, x>0\}.  
\end{equation*}
The well-posedness in $H_{\eq}^1(\Gamma)$ follows analogously to  \cite[Lemma 2.6]{Gol19}.
We use index $\cdot_{\eq}$  to   denote counterparts of the objects for the space $H_{\eq}^1(\Gamma)$. 

It is known that $d^0_{\omega,\eq}=S^0_{\omega}(\phi_\gamma)$ (see page 12 in \cite{NataMasa2020}) for any $\gamma\in\mathbb{R}$, where 
$$\phi_\gamma(x)=\left(\left\{\tfrac{(p+1)\omega}{2}\mathrm{sech}^2\left(\tfrac{(p-1)\sqrt{\omega}}{2}x+\mathrm{arctanh}(\tfrac{\gamma}{N\sqrt{\omega}})\right)\right\}^{\frac{1}{p-1}}\right)_{e=1}^N.$$
Then for $0<\omega_{0,\eq}<\omega$ (observe that $\omega_{0,\eq}\leq \omega_0$) one can repeat all the proofs in Section 3 and  Subsections \ref{main_res} and \ref{proof_cor} with $H_{\eq}^1(\Gamma)$ instead of $H^1(\Gamma).$ Thus, we get the following results.

\begin{proposition}
Let $p>1, \gamma\in \mathbb{R},$  $\omega>\omega_{0,\eq}$. If  $V(x)=\overline{V(x)}$ is symmetric and satisfies Assumptions 1-3, then the set of ground states $\mathcal{G}_{\omega,\eq}$ is not empty, in particular,  $\mathcal{G}_{\omega,\eq}=\mathcal{M}_{\omega,\eq}$. If $\varphi_{\omega}\in\mathcal{G}_{\omega,\eq}$, then there exist $\theta\in\mathbb{R}$ and a positive function $\phi\in H^{1}_{\eq}(\Gamma)$ such that $\varphi_{\omega}(x)=e^{i\theta}\phi(x)$.
\end{proposition}

\begin{theorem} 
Let $p>5$,  $\gamma\in \mathbb{R}$,  $\omega>\omega_{0,\eq}$.  If  $V(x)=\overline{V(x)}$ is symmetric and satisfies Assumptions 1-4, $\varphi_{\omega}(x)\in\mathcal{G}_{\omega,\eq}$, and  $\partial^{2}_{\lambda}E(\varphi^{\lambda}_{\omega})\left|_{\lambda=1}\right.<0$, then the standing wave solution $e^{i\omega t}\varphi_{\omega}(x)$ of \eqref{NLS} is orbitally unstable in $H^{1}_{\eq}(\Gamma)$ and therefore in $H^{1}(\Gamma)$.
\end{theorem}

\begin{corollary}
Assume that $V(x)=\dfrac{-\beta}{x^\alpha}$,  $\beta>0$,  $0<\alpha<1$, $\gamma\in \mathbb{R}$. Let $p>5$ and $\varphi_{\omega}(x)\in\mathcal{G}_{\omega,\eq}$. Then there exists $\omega^*_{\eq}\in(\omega_{0,\eq},\infty)$ such that for any $\omega\in (\omega^*_{\eq},\infty)$ the standing wave solution $e^{i\omega t}\varphi_{\omega}(x)$ of \eqref{NLS} is orbitally unstable in $H^{1}(\Gamma)$.  
\end{corollary}
\begin{remark}
$(i)$\, Observe that when dealing with $H^1_{\eq}(\Gamma)$, no restriction on $\gamma$ appears. This is due to the fact that  the corresponding constrained variational problem is closely related to  the one on $\mathbb{R}$, which in turn admits a minimizer for any $\gamma$ (see \cite[Remark 3.1]{NataMasa2020}).

$(ii)$\,
Consider 
 $$i\partial_{t}u=-\partial_x^2u-\gamma\delta(x)u+V(x){u}-\left|u\right|^{p-1}u,\quad (t,x)\in \mathbb{R}\times\mathbb{R},$$
 $\gamma\in \mathbb{R}$. Notice that the above results are valid with $H^{1}_{\eq}(\Gamma)$ substituted by\newline $H^{1}_{\mathrm{rad}}(\mathbb{R})=\{f\in H^1(\mathbb{R}): f(x)=f(-x)\}$ and analogous assumptions on $V(x)$. One only needs to recall that $d^0_{\omega,\rad}=S^0_{\omega}(\phi_\gamma)$ (see \cite[Theorem 1]{RJJ}), where \newline
 $\phi_\gamma(x)=\left\{\frac{(p+1)\omega}{2}\mathrm{sech}^2\left(\frac{(p-1)\sqrt{\omega}}{2}|x|+\mathrm{arctanh}(\frac{\gamma}{2\sqrt{\omega}})\right)\right\}^{\frac{1}{p-1}}$.
\end{remark}


\section*{Appendix}

Below we show  some properties of the operator $ H_{\gamma, V}$  introduced by \eqref{Ku}.

\begin{lemma} \label{Ape1}
Let $\gamma\in \mathbb{R}$ and $V(x)=\overline{V(x)}\in L^{1}(\Gamma)+L^{\infty}(\Gamma)$. The quadratic form $ {F}_{\gamma, V}$  given by \eqref{F} is semibounded and closed, and the operator $ H_{\gamma, V}$ defined by
\begin{equation*}
\begin{split}
&( H_{\gamma, V}v)_{e}=-v''_{e}+V_ev_{e},\\
&\dom( H_{\gamma, V})=\left\{v\in H^{1}(\Gamma):\,\,\,-v''_{e}+V_ev_{e}\in L^{2}(\mathbb{R}^{+}),\,\,\,\sum^{N}_{e=1}v^{\prime}_{e}(0)=-\gamma v_{1}(0)\right\}.
\end{split}
\end{equation*}
 is the self-adjoint operator associated with $ {F}_{\gamma, V}$  in $L^{2}(\Gamma)$. 
\end{lemma}
\begin{proof}
We can write $V(x)=V_{1}(x)+V_{2}(x)$, with $V_{1}\in L^{1}(\Gamma)$ and $V_{2}\in L^{\infty}(\Gamma)$. Thus, using the Gagliardo-Nirenberg inequality (see formula (2.1) in \cite{CAFiNo2017}) and the Young inequality, we have
\begin{equation}\label{ee2}
\begin{split}
\left|\int\limits_{\Gamma}V(x)\left|v(x)\right|^{2}dx\right| &\leq\left\|V_{1}\right\|_1\left\|v\right\|^{2}_\infty+\left\|V_{2}\right\|_\infty\left\|v\right\|^2_2\\
&\leq C\left\|V_{1}\right\|_1\left\|v'\right\|_2\left\|v\right\|_2+\left\|V_{2}\right\|_\infty\left\|v\right\|^2_2\\
&\leq \epsilon\left\|v'\right\|^2_2+C_{\varepsilon}\left\|v\right\|^2_2,\,\,\,\epsilon>0.
\end{split}
\end{equation}
Similarly, by  the Sobolev embedding, we obtain
\begin{align*}\label{ee1}
\left|\gamma\left|v_{1}(0)\right|^{2}\right|\leq |\gamma|\|v\|_\infty^2\leq C\|v'\|_2\left\|v\right\|_2\leq\epsilon\left\|v'\right\|^2_2+C_{\epsilon}\left\|v\right\|^2_2.
\end{align*}
Therefore,
\begin{equation}\label{l_4.9.5}
\left|\gamma\left|v_{1}(0)\right|^{2}+\int\limits_{\Gamma}V(x)|v(x)|^2dx\right|\leq2\epsilon\left\|v'\right\|^2_2+C_\varepsilon\left\|v\right\|^2_2,\,\,\mbox{for every}\,\,\epsilon>0.
\end{equation}
Then, by the KLMN theorem \cite[Theorem X.17]{ReeSim75}, we infer that the quadratic form   $ {F}_{\gamma, V}$ is associated with a semibounded self-adjoint operator $ T_{\gamma, V}$ defined by (observe that $A=H_{0,0}$ in \cite[Theorem X.17]{ReeSim75}, i.e. $V\equiv 0, \gamma=0$)
\begin{equation*}
\begin{split}
&\dom( T_{\gamma, V})=\left\{u\in H^1(\Gamma):\,\exists\, y\in L^2(\Gamma)\, s.t.\, \forall v\in H^1(\Gamma),\,\,   F_{\gamma, V}(u, v)=(y, v)_2\right\},\\
&  T_{\gamma, V}u=y.
\end{split}    
\end{equation*}
It is easily seen that $\dom( H_{\gamma, V})\subseteq\dom( T_{\gamma, V})$ and $ T_{\gamma, V}u= H_{\gamma, V}u, \,\,u\in\dom( H_{\gamma, V})$. Hence it is sufficient to prove that $\dom( T_{\gamma, V})\subseteq\dom( H_{\gamma, V})$.

Let $\tilde u\in \dom( T_{\gamma, V})$ and $\tilde v\in H^1(\Gamma)$, then there exists $\tilde y\in L^2(\Gamma)$ such that 
\begin{equation}\label{l_4.9.1}
F_{\gamma, V}(\tilde u, \tilde v)=\int\limits_\Gamma (\tilde u' \overline{\tilde v'}+V\tilde u\overline{\tilde v})dx-\gamma \tilde u_1(0)\overline{\tilde v_1(0)}= (\tilde y, \tilde v)_2.  
\end{equation}
Observe that $\tilde y-V\tilde u\in L^1_{loc}(\Gamma)$ and set $$z=(z_e)_{e=1}^N,\quad z_e(x)=\int\limits_0^x\left(\tilde y_e(t)-V_e(t)\tilde u_e(t)\right)dt.$$ 

Suppose now additionally that $\tilde v$ has a compact support, then 
\begin{equation}\label{l_4.9.2}
\int\limits_\Gamma(\tilde y-V\tilde u)\overline{\tilde v}dx=\int\limits_\Gamma z'\overline{\tilde v}dx=-\overline{\tilde v_1(0)}\sum\limits_{e=1}^Nz_e(0)-\int\limits_\Gamma z \overline{\tilde v'}dx.   
\end{equation}
From \eqref{l_4.9.1} we deduce
\begin{equation}\label{l_4.9.3}
  \int\limits_\Gamma(\tilde y-V\tilde u)\overline{\tilde v}dx=\int\limits_\Gamma \tilde u' \overline{\tilde v'}dx-\gamma \tilde u_1(0)\overline{\tilde v_1(0)}. 
\end{equation}
Combining \eqref{l_4.9.2} and \eqref{l_4.9.3} we get 
\begin{equation}\label{l_4.9.4}
\int\limits_\Gamma(\tilde u'+z)\overline{\tilde v'}dx+\overline{\tilde v_1(0)}\left(-\gamma\tilde u_1(0)+\sum\limits_{e=1}^Nz_e(0)\right)=0.    
\end{equation}
Choose $\tilde v=(\tilde v_e)_{e=1}^N$ such that $\tilde v_1(x)\in C_0^\infty (\mathbb{R}^+)$ and $\tilde v_2(x)\equiv\ldots\equiv\tilde v_N(x)\equiv 0.$ 
Then we obtain
$$\int\limits_0^\infty(\tilde u'_1+z_1)\overline{\tilde v'_1}dx=0,$$ therefore $\tilde u'_1+z_1\equiv const\equiv c_1$. We have used that $\tilde {u}'_1+z_1\in\mathrm{Ran}(A)^{\perp},$ where $Av=v'$ with $\dom(A)=C_0^\infty(\mathbb{R}^+)$ in $L^2(\mathbb{R}^+)$. Analogously $\tilde u'_e+z_e\equiv const\equiv c_e,\,\, e=2,\ldots, N. $
Finally, from \eqref{l_4.9.4} we deduce
$$\overline{\tilde v_1(0)}\left(-\gamma\tilde u_1(0)-\sum\limits_{e=1}^N(\tilde u'_e(0)+z_e(0))+\sum\limits_{e=1}^Nz_e(0)\right)=0. $$
Assuming that $\tilde v_1(0)\neq 0,$ we arrive at $\sum\limits_{e=1}^N\tilde u'_e(0)=-\gamma \tilde u_1(0).$ Moreover, $-\tilde u''+V\tilde u=z'+V\tilde u=\tilde y-V\tilde u+V\tilde u=\tilde y\in L^2(\Gamma).$ Hence $\tilde u\in \dom( H_{\gamma, V})$ and $\dom( T_{\gamma, V})\subseteq \dom(H_{\gamma, V}).$
\end{proof}
\begin{lemma}\label{H^2} Suppose that $V(x)=\overline{V(x)}\in L^2_{\varepsilon}(\Gamma)+L^\infty(\Gamma)$, i.e. for any $\varepsilon>0$ and $V\in L^2_{\varepsilon}(\Gamma)+L^\infty(\Gamma)$ there exists a representation $V=V_1+V_2, \, V_1\in L^2(\Gamma), V_2\in L^{\infty}(\Gamma)$, with $\|V_1\|_2^2\leq \varepsilon$.   Then we have
\begin{equation} \label{D_2}\dom( H_{\gamma, V})=\left\{v\in H^{1}(\Gamma):\,\,\,v_{e}\in H^{2}(\mathbb{R}^{+}),\,\,\,\sum^{N}_{e=1}v^{\prime}_{e}(0)=-\gamma v_{1}(0)\right\}:=D_{H^2}.
\end{equation}
Moreover,  for $m$ sufficiently large,  $ H_{\gamma, V}$-norm  $\|(H_{\gamma, V}+m)\cdot\|_2$ is equivalent to $H^2$-norm on $\Gamma.$
\end{lemma}
\begin{proof}
Observe that, by  $ V(x)\in L^2_{\varepsilon}(\Gamma)+L^\infty(\Gamma)$,  the Sobolev and the Young inequalities we get
\begin{equation}\label{r.4.10.2}
 \|Vv\|_2^2\leq  \|V_1\|_2^2\|v\|_\infty^2+\|V_2\|_\infty^2\|v\|_2^2\leq\varepsilon\|v\|_{H^2(\Gamma)}^2+C\|v\|_2^2
\end{equation}
and, 
\begin{equation}\label{r.4.10.1}
\begin{split}
  &|(v'', Vv)_2|\leq \|v''\|_2  \|Vv\|_2\leq \|v''\|_2 \|V_1\|_2\|v\|_\infty+  \|v''\|_2 \|V_2\|_\infty\|v\|_2\\&\leq 
  C_1\|v''\|_2\|V_1\|_2\|v\|_{H^2(\Gamma)}+C_2\|v''\|_2 \|v\|_2 \leq \varepsilon \|v\|^2_{H^2(\Gamma)}+\varepsilon \|v''\|^2_2+C_\varepsilon \|v\|_2^2\\&
   \leq 2\varepsilon \|v\|^2_{H^2(\Gamma)}+C_\varepsilon \|v\|_2^2.
  \end{split}
\end{equation}
It is immediate from \eqref{r.4.10.2}, \eqref{r.4.10.1} that 
\begin{equation}\label{H21} \| H_{\gamma, V}v\|_2^2=\|v''\|_2^2+2\re(v'', Vv)_2+ \|Vv\|_2^2\leq C_1\|v\|^2_{H^2(\Gamma)}.\end{equation}
And for $m$ sufficiently large, inequalities   \eqref{r.4.10.2} and \eqref{r.4.10.1} imply
\begin{equation}\label{H22}\| H_{\gamma, V}v\|_2^2+m^2\|v\|_2^2=\|v''\|_2^2+2\re(v'', Vv)_2+ \|Vv\|_2^2+m^2\|v\|_2^2\geq C_2\|v\|^2_{H^2(\Gamma)}.\end{equation}
Thus, we get \eqref{D_2}.

The second assertion follows from \eqref{H21},\eqref{H22}, and 
\begin{align*}
& \|( H_{\gamma, V} + m)v\|_2^2=\| H_{\gamma, V}v\|_2^2+m^2\|v\|_2^2+2m(H_{\gamma, V}v,v)_2,\quad   \\
& |(H_{\gamma, V}v,v)_2|\leq  \|H_{\gamma, V}v\|_2\|v\|_2\leq \varepsilon \|H_{\gamma, V}v\|_2^2+C_\varepsilon\|v\|_2^2.
\end{align*}
\end{proof}
\begin{remark}
Observe that there exists potential $V(x)$ satisfying Assumptions 1-4 such that $\dom( H_{\gamma, V})\neq D_{H^2}.$  For  example, consider $V(x)=-1/x^\alpha,\,1/2\leq \alpha<1,$ and $N=\gamma=2$, then   $v=(e^{-x}, e^{-x})\in D_{H^2}$, but 
\begin{align*}
 \|H_{\gamma, V}v\|_2^2=2\|-v''_1-\frac{v_1}{x^\alpha}\|_2^2>2e^{-2\varepsilon}\int\limits_0^\varepsilon\frac{dx}{x^{2\alpha}}=\infty.
\end{align*}
\end{remark}
\begin{lemma}\label{lemma_4.11} Let $\gamma>0$ and $V(x)=\overline{V(x)}$ satisfy Assumptions 1 e 3. Then the following assertions hold.

\noindent$(i)$ The number $-\omega_0$ defined by \eqref{omega_0} is negative.

\noindent$(ii)$\, Let also $m>\omega_0$, then  $\sqrt{ {F}_{\gamma, V}(v)+m\|v\|_2^2}$ defines a norm equivalent to the $H^1$-norm.

\noindent$(iii)$\, The number $-\omega_0$ is the first eigenvalue of $H_{\gamma, V}$. Moreover, it is simple,   and there exists the corresponding  positive eigenfunction $\psi_0\in\dom(H_{\gamma, V})$, i.e. $H_{\gamma, V}\psi_0=-\omega_0\psi_0.$ 
\end{lemma}
\begin{proof}
$(i)\,$
To show $-\omega_0<0$, observe that 
 \begin{align}\label{F_gamma}
-\omega_0=\inf\sigma(H_{\gamma,V})=\inf\left\{ {F}_{\gamma, V}(v)\,:\,v\in H^{1}(\Gamma),\,\,\left\|v\right\|^2_2=1\right\}.
\end{align} Consider 
 $v^{\lambda}(x)=\lambda^{\frac{1}{2}}v(\lambda x)$ with $\lambda>0$. Hence
\begin{align*}
 {F}_{\gamma, V}(v^{\lambda})&=\lambda^{2}\left\|v'\right\|^2_2-\lambda\gamma\left|v_{1}(0)\right|^{2}+(Vv^{\lambda}, v^{\lambda})_2.
\end{align*} 
For $\lambda$ small enough, we have $ {F}_{\gamma, V}(v^{\lambda})<0$. Finally, $-\omega_0$ is finite since   ${F}_{\gamma, V}(v)$ is lower semibounded.

$(ii)\,$  Let $\varepsilon>0$. Firstly, notice  that from \eqref{l_4.9.5} one easily gets
$$ {F}_{\gamma, V}(v)+m\|v\|_2^2\leq (1+2\varepsilon)\|v'\|_2^2+(C_\varepsilon+m)\|v\|_2^2\leq C_1\|v\|^2_{H^1(\Gamma)}.$$
Secondly, for $\varepsilon$ and $\delta$ sufficiently small, 
\begin{align*}
 & {F}_{\gamma, V}(v)+m\|v\|_2^2= \delta\|v'\|_2^2+(1-\delta)\left(\|v'\|_2^2+\frac{1}{1-\delta}(Vv,v)_2-\frac{\gamma}{1-\delta}|v_1(0)|^2\right)+m\|v\|_2^2\\&\geq\delta\|v'\|_2^2-(1+\varepsilon)(1-\delta)\omega_0\|v\|_2^2+m\|v\|_2^2\geq C_2\|v\|^2_{H^1(\Gamma)}.   
\end{align*}
Indeed, the family of sesquilinear forms $$\mathrm{t}(\kappa)[u,v]=(u',v')_2+\frac{1}{1-\kappa}(Vu,v)_2-\frac{\gamma}{1-\kappa}(u_1(0)\overline{v_1}(0))$$ is holomorphic of type (a) in the sense of Kato in the complex neighborhood of zero (see \cite[Chapter VII, \S 4]{Kat95} for the definition and \cite[Chapter VI, \S 1, Example 1.7]{Kat95} for the proof of sectoriality). Using inequality $(4.7)$ in \cite[Chapter VII]{Kat95}  with $\kappa=\kappa_2=0, \kappa_1=\delta$, we obtain $|\mathrm{t}(\delta)[v]-\mathrm{t}(0)[v]|\leq \varepsilon|\mathrm{t}(0)[v]|.$ Hence
$$\mathrm{t}(\delta)[v]\geq \mathrm{t}(0)[v]-\varepsilon|\mathrm{t}(0)[v]|= {F}_{\gamma, V}(v)-\varepsilon| {F}_{\gamma, V}(v)|\geq-(1+\varepsilon)\omega_0\|v\|_2^2.$$
$(iii)\,$  \textit{Step 1.}   
Let $\{v_{n}\}$ be a minimizing sequence, that is, ${F}_{\gamma, V}(v_{n})\underset{n\to\infty}{\longrightarrow} -\omega_{0}$, $\left\|v_{n}\right\|^{2}_2=1$ for all $n\in\mathbb{N}$. From $(ii)$, we deduce that $\{v_{n}\}$ is bounded in $H^{1}(\Gamma)$. Then there exist a subsequence $\{v_{n_k}\}$ of $\{v_n\}$ and $v_0\in H^{1}(\Gamma)$ such that $\{v_{n_k}\}$ converges weakly to $v_0$ in $H^1(\Gamma)$. Observe that, by the weak lower semicontinuity of $L^2$-norm and ${F}_{\gamma, V}(\cdot)$, we get  $\left\|v_0\right\|_2\leq1$ and   
\begin{align*}
F_{\gamma, V}(v_0)\leq\lim_{k\rightarrow\infty}F_{\gamma, V}(v_{n_k})=-\omega_{0}<0.
\end{align*}
We have $\left\|v_0\right\|_{2}=1$, since, otherwise, there would exist $\lambda>1$ such that $\left\|\lambda v_0\right\|_{2}=1$ and $F_{\gamma, V}(\lambda v_0)=\lambda^2F_{\gamma, V}(v_0)<-\omega_0$, which is a contradiction.  Consequently  $v_0$ is a minimizer for \eqref{F_gamma}. 

Let $\psi_{0}=\left|v_0\right|$, then $\psi_{0}\geq0$ on $\Gamma$  and $\left\|\psi_{0}\right\|^{2}_2=\left\|v_{0}\right\|^{2}_2=1$. Notice that
$
\left\|\psi'_{0}\right\|^{2}_2\leq\left\|v'_0\right\|^{2}_2,$
therefore ${F}_{\gamma, V}(\psi_{0})\leq{F}_{\gamma, V}(v_0)$.  Then $\psi_{0}$ is a minimizer of \eqref{F_gamma}.
This implies the existence of  the Lagrange multiplier $-\mu$ such that  
$$F'_{\gamma, V}(\psi_0)=-\mu Q'(\psi_0),\quad Q(v)=\|v\|_2^2.$$ Repeating  the arguments from the proof of \cite[Theorem 4]{AQFF}, we get $\psi_0\in \dom(H_{\gamma, V})$ and 
$$H_{\gamma, V}\psi_0=-\mu \psi_0.$$
Multiplying the above equation by $\overline{\psi_0}$ and integrating we conclude $\mu=\omega_0.$  Recalling  that  $V(x)\leq 0$ a.e. on $\Gamma$, and  arguing as in the proof of Proposition \ref{GS}, one can show that $\psi_0>0$ on $\Gamma.$ Notice that one needs to apply  \cite[Theorem 1]{LVL} with $\beta(s)=\omega_0 s.$

\textit{Step 2.} Suppose that $u_0$ is a  nonnegative solution of
\begin{equation}\label{eigen} 
H_{\gamma, V}u_0=-\omega_0 u_0.
\end{equation}Let us show that there exists $C>0$ such that $u_0(x)=C\psi_0(x)$. Assume that this is false. Then there exists $C>0$ such that $\widetilde{u}_0(x)=u_0(x)-C\psi_0(x)$ takes both positive and negative values. We have $H_{\gamma, V}\widetilde{u}_0=-\omega_0 \widetilde{u}_0,$ consequently $\widetilde{v}_0=\widetilde{u}_0/\|\widetilde{u}_0\|_2$ is the minimizer of \eqref{F_gamma}. Arguing as in \textit{Step 1}, one can show that $|\widetilde{v}_0|$ is also a  minimizer and $|\widetilde{v}_0|>0$. Therefore, $\widetilde{u}_0(x)$ has a constant sign. This is a contradiction.   

Suppose now that $u_0$ is an arbitrary solution to \eqref{eigen} such that $\|u_0\|_2^2=1$ (that is $u_0$ is a minimizer of \eqref{F_gamma}).  Define $w_0=|\re u_0|+i|\im u_0|,$ then $|w_0|=|u_0|$ and $|w'_0|=|u'_0|$, consequently $F_{\gamma, V}(u_0)=F_{\gamma, V}(w_0)$
and $\|w_0\|_2^2=1.$ Therefore, $w_0$ is  a minimizer of \eqref{F_gamma}. This implies that $w_0$ satisfies \eqref{eigen}, and, in particular,  $|\re u_0|$ and $|\im u_0|$ satisfy \eqref{eigen}. Thus, $|\re u_0|=C_1 \psi_0$ and $|\im u_0|=C_2\psi_0,\, C_1,C_2>0$, consequently $\re u_0=\widetilde{C}_1 \psi_0$ and $\im u_0=\widetilde{C}_2\psi_0,\, \widetilde{C}_1,\widetilde{C}_2\in\mathbb{R},$ since $\re u_0$ and $\im u_0$ do not change the sign. Finally, $u_0=\widetilde{C}_1\psi_0+i\widetilde{C}_2\psi_0=\widetilde{C}\psi_0, \, \widetilde{C}\in \mathbb{C},$ and therefore ,$-\omega_0$ is simple.

\end{proof}

\noindent 
{\bf Acknowledgments.} The authors are kindly grateful to Prof. Gláucio Terra for  the proof of Remark \ref{V_neg}.

\end{document}